\documentclass[12pt]{article}

\usepackage[utf8]{inputenc}
\usepackage[T1]{fontenc}
\usepackage[margin=1in]{geometry}
\usepackage{setspace}
\usepackage{natbib}
\usepackage{hyperref}

\usepackage{xcolor}
\usepackage{graphicx}
\usepackage{caption}
\usepackage{subcaption}
\usepackage{float}
\usepackage{booktabs}
\usepackage{algorithm}
\usepackage{algorithmicx}
\usepackage{algpseudocode}
\usepackage{tikz-cd}
\usepackage{wrapfig}
\usepackage{authblk}
\usepackage{enumerate}
\usepackage{amsmath, amssymb, amsfonts, amsthm, mathtools}
\usepackage{mathrsfs}
\usepackage[title]{appendix}
\usepackage{listings}
\usepackage{comment}
\usepackage{scalerel}
\usepackage{chngcntr}
\usepackage{apptools}
\usepackage[]{lineno}
\usepackage{float} % in the preamble
\usepackage{placeins} % preamble
\usepackage{graphicx}
\graphicspath{{fig/}{fig/2-Datasets/}{fig/}{fig/}{fig/2-Datasets/}}

\makeatletter
\newtheorem*{rep@theorem}{\rep@title}
\newcommand{\newreptheorem}[2]{%
\newenvironment{rep#1}[1]{%
\def\rep@title{#2 \ref{##1}}%
\begin{rep@theorem}}%
{\end{rep@theorem}}}
\makeatother

\newreptheorem{theorem}{Theorem}

\theoremstyle{definition}

\AtAppendix{\counterwithin{theorem}{section}}

\definecolor{darkblue}{rgb}{0.0, 0.0, 0.8}
\definecolor{darkred}{rgb}{0.8, 0.0, 0.0}
\definecolor{darkgreen}{rgb}{0.0, 0.5, 0.0}
\definecolor{ncolor}{rgb}{0.8, 0.8, 0.0}

\newcommand{\keywords}[1]{%
  \vspace{1em}
  \noindent\textbf{Keywords: }#1
}

\title{Improving Topological Detection of Weather Regimes in Climate Dynamical Systems}

\author{
Soheil Anbouhi\\
Department of Mathematics and Computer Science\\
Western Carolina University\\
Cullowhee, NC 28723, USA\\
\texttt{sanbouhi@wcu.edu}
}

\date{\today}

\begin{document}

\maketitle

%========================================
% Abstract
%========================================
\begin{abstract}

Weather regimes provide a useful framework for describing large-scale atmospheric variability and its impacts on regional weather. Despite extensive study, there is no universally accepted definition or method for identifying weather regimes. Recent work has shown that regimes can be interpreted geometrically as topological structures in atmospheric phase space, allowing their detection without prescribing the number of regimes in advance. In this framework, regimes are identified using a density--radius bifiltration combined with persistent homology, a well-established tool in topological data analysis. A limitation of this approach is its reliance on density estimation, which can over-smooth weakly populated yet dynamically meaningful regions of phase space. This limitation is particularly relevant for datasets with weakly separated or sparsely sampled regimes, such as the southern jet regime in the North Atlantic and the thin zonal loops of the Charney--DeVore system. Here, we introduce a centrality--radius bifiltration that emphasizes local connectivity in phase space and overcomes the limitations of the density-based method, in particular by recovering the southern jet regime and the thin loop structures of the Charney--DeVore system. The resulting regime-related topological structures are robust across a range of analysis scales, suggesting that they reflect the intrinsic organization of atmospheric phase space rather than artifacts of parameter choice. The method provides a practical and reproducible diagnostic tool for identifying weather regimes across datasets. We also integrate our method into an existing analysis framework.

\end{abstract}

%==========================
%Keywords
%==========================

%==========================
%Introduction
%==========================
\section{Introduction}\label{sec:intro}

The atmosphere is a nonlinear, chaotic dynamical system. Yet, its large-scale variability can be described using a few recurrent, quasi-stationary circulation patterns, known as weather regimes. These regimes reflect the atmosphere’s tendency to remain in particular circulation states for extended periods before transitioning to another \citep{legras1985persistent, vautard1990multiple, corti1999signature, hannachi2017weather}.

In phase space—where each point represents a possible atmospheric state—weather regimes correspond to regions in which the system’s trajectory slows, lingers for extended periods (quasi-stationary), and to which it returns repeatedly (recurrent). In the Northern Hemisphere, a familiar case is atmospheric blocking, where a stable high-pressure anomaly disrupts the typical zonal flow for several days. Similarly, large-scale patterns such as the North Atlantic Oscillation (NAO) and the Pacific–North American (PNA) pattern describe recurrent shifts in pressure and circulation that shape weather across broad regions \citet{hannachi2017low}.

The importance of weather regimes lies in their ability to simplify the complexity of atmospheric behavior while improving predictability by linking large-scale circulation patterns to local weather impacts \citet{hannachi2017low}. Despite decades of research, however, there is still no universally accepted definition or detection method for weather regimes. Traditional approaches—such as clustering techniques, Markov-based models, or self-organizing maps—rely on user-defined parameters, often including or implicitly determining the number of regimes. As a result, applications to Euro-Atlantic circulation have reported between two and seven distinct regimes depending on the method and parameter choices \citet{strommen2023topological}. Moreover, many of these approaches become computationally inefficient as the dimensionality of the data increases, reflecting the well-known curse of dimensionality.

Motivated by these challenges, \citet{strommen2023topological} (hereafter S23) proposed a topology-based framework for identifying weather regimes. In this framework, regimes correspond to non-trivial topological structures—such as connected components and loops—within the attractor of the atmospheric system. To detect these structures, the authors employ persistent homology (PH), a key tool in topological data analysis \citet{carlsson2009topology}.

PH tracks topological features across spatial scales and quantifies their persistence by recording the scales at which they appear and disappear. This information is summarized in a diagram, where each feature is represented by its birth and death values. Features with longer persistence are typically interpreted as robust structures, while short-lived features are more sensitive to the choice of scale and may reflect noise or fine-scale variability.

Because climate datasets are typically high-dimensional and densely sampled, some meaningful topological features may disappear when standard PH is applied. To address this issue, S23 introduced a density–radius bifiltration in which the data is first filtered by a density function (e.g., Gaussian kernel density estimation), then PH is applied to extract topological features. This allows one to track how these features evolve as the density parameter varies.

This method was tested on three classical dynamical systems and an observational climate dataset that exhibit very different regime behavior. These dynamical systems differ in how dense and how persistent their regimes are, making them a strong test of how the method generalizes. In all cases, the extracted topological features successfully describe the essential regime behavior. Therefore, topology can be considered as a unifying property across these most well-known regime systems. Indeed, S23, argued that no single, simple definition of a regime—based solely on density, persistence, or recurrence—can adequately describe all systems.

Furthermore, their method overcomes several limitations of earlier methods: it does not require the number of regimes to be set in advance and scales efficiently to large, high-dimensional datasets.

A limitation of their method is that, in some cases, the choice of density estimator can significantly influence regime identification. In S23, two estimators are considered: Gaussian kernel density estimation (KDE) and direct binning method (histogramming). KDE can be viewed as a smooth, continuous analogue of direct binning. Both methods require careful tuning of parameters (the bandwidth in KDE or the number of bins in direct binning) to avoid under-smoothing noise or over-smoothing meaningful structure. Although automatic bandwidth-selection rules such as Scott’s rule exist for KDE, the resulting Gaussian KDE–based bifiltration may still over--smooth small-scale yet dynamically significant features. This limitation is evident when studying the North Atlantic jet latitude index, which exhibits a trimodal distribution corresponding to southern, central, and northern jet positions. The southern jet regime, however, is not well separated \citep{woollings2010variability, strommen2023topological} and therefore fails to appear as a distinct connected component in the KDE-based bifiltration (see Section~\ref{subsec:jetlat}).

A similar issue arises in the Charney–DeVore system. Its attractor exhibits two regimes, one of which contains a set of low-density loops that spiral outward in phase space and are clearly visible by eye (see Fig.~\ref{fig:CDV_comparison}). The Gaussian KDE tends to smooth over these thin structures and, as a result, the topological method fails to detect them. While the direct binning approach improves the detection of these features, it remains less reliable, being highly sensitive to bin size, bin alignment, and data dimensionality \citep{Scott1992, Silverman1986}.

To overcome these limitations, we replace density with a measure of centrality in the bifiltration framework. Instead of estimating density directly, we measure how strongly each state is connected to its nearby neighbors. This is achieved using \textit{distance-to-measure} functions, which compute the average distance to a point’s $k$ nearest neighbors and are known to be robust to outliers \citet{chazal2011geometric}. The parameter $k$ controls the scale of locality: smaller values emphasize fine-scale structure, while larger values provide increased robustness to noise.

The proposed centrality-based bifiltration improves upon the density–radius approach for the North Atlantic jet latitude dataset by successfully capturing all three jet regimes (southern, central, and northern) across a wide range of k values. It also enhances the detection of thin, low-density loop structures associated with the zonal regime in the Charney–DeVore system, which are not reliably identified by the KDE-based method. In addition, the method is computationally more efficient than the KDE-based bifiltration and, like the KDE-based approach, scales to large, high-dimensional datasets without requiring the number of regimes to be specified in advance.

%=======================Outline=========================
\subsection*{Outline}\label{subsec:outline}
%====================================================

Section~\ref{sec:data} introduces the datasets analyzed in this study and briefly reviews their associated weather regimes. Section~\ref{sec:method} provides background on persistent homology, presents our proposed method, describes the computational setup, and discusses the significance tests used to evaluate the robustness of the results. Section~\ref{sec:results} presents the results and compares them with those of S23. Finally, Section~\ref{sec:discussion} summarizes the main findings, discusses the strengths and limitations of the approach, and outlines directions for future work.

%==========================
%Datasets
%==========================
\section{Datasets}\label{sec:data}
We use the same datasets analyzed by S23: three from idealized dynamical models and one from reanalysis data. We begin with Lorenz--63 as a reference system before focusing on the Charney--DeVore system and the North Atlantic jet latitude dataset. Additionally, results for the Lorenz--96 system, which are included for comparison with S23, are provided in the Appendix. Figure~\ref{fig:datasets-all} shows the datasets considered in this study. For full details, see Appendix B of S23.

\begin{figure}[ht]
\centering
\begin{subfigure}{0.29\linewidth}
    \centering
    \includegraphics[width=\linewidth]{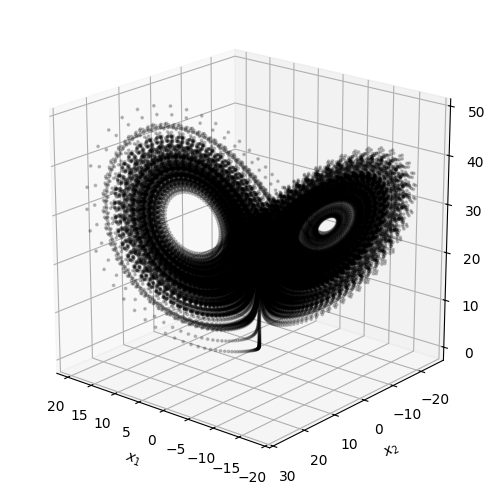}
    \caption{\small Lorenz--63 attractor}
    \label{fig:dataset-L63}
\end{subfigure}
\hspace{0.3cm}
\begin{subfigure}{0.30\linewidth}
    \centering
    \includegraphics[width=\linewidth]{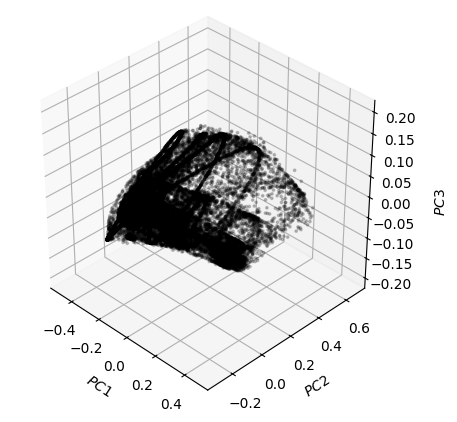}
    \caption{\small Charney--DeVore attractor}
    \label{fig:dataset-CdV}
\end{subfigure}
\hspace{0.3cm}
\begin{subfigure}{0.3\linewidth}
    \centering
    \includegraphics[width=\linewidth]{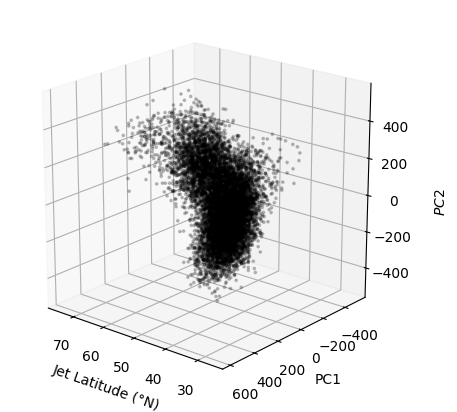}
    \caption{\small North Atlantic jet dataset}
    \label{fig:dataset-JetLat}
\end{subfigure}

\caption{\small Datasets used in this study: (a) Lorenz--63 attractor with 40{,}000 data points, (b) Charney--DeVore attractor projected onto the first three principal components with 40{,}000 data points, and (c) North Atlantic jet latitude dataset projected onto three dimensions with approximately 9{,}800 data points.}
\label{fig:datasets-all}
\end{figure}

%==========================L63==========================
\subsection{Lorenz--63 model}
%====================================================

The Lorenz–63 system \citet{lorenz1963deterministic} is a three-variable chaotic model originally developed as a simplified representation of atmospheric convection. Its attractor has the well-known ``butterfly'' shape, with two spiral lobes often interpreted as alternating regimes analogous to zonal and blocked flow patterns in the atmosphere \citep{palmer1993extended,baines2008lorenz,molteni2019heuristic,strommen2023topological} (see Fig.~\ref{fig:datasets-all}(a) ).

%==========================CdV==========================
\subsection{Charney--DeVore model}
%====================================================
The Charney--DeVore (CdV) model is a simple six-variable system used to study large-scale atmospheric flow \citet{charney1979multiple}. It produces two regimes corresponding to blocked and zonal flow. In the classical view \citet{charney1979multiple}, these regimes are associated with two fixed points around which the model’s trajectory slows down. The blocked regime is more stable and quasi-stationary, while the zonal regime is more turbulent and chaotic \citet{pomeau1980intermittent}. In phase space, the fixed point for the blocked regime lies near the dense central region of the attractor, while the fixed point for the zonal regime is found near the back left corner, and appears as a set of low-density loops spiraling outward across a wider area \citet{strommen2023topological}. These loops can be interpreted as weak bursting transients and are visible in the raw data. Fig.~\ref{fig:datasets-all}(b) shows the projection of this dataset onto the first three principal components.

%==========================Jet==========================
\subsection{North Atlantic jet}
%====================================================

The final dataset is observational. We use a three-dimensional representation of North Atlantic winter jet variability, which we refer to as \texttt{JetLat}, following S23. The dataset is derived from the ERA20C reanalysis \citet{poli2016era} (1900--2010) and validated against ERA-Interim \citet{dee2011era} (1979--2015). For each winter day (December--February), it includes the daily jet latitude index—the latitude of maximum 850~hPa zonal wind—and the first two principal components of 850~hPa zonal wind anomalies (see Fig.~\ref{fig:datasets-all}(c)). Details on how the jet latitude is computed are given in \citep{parker2019seasonal, strommen2020jet}.

The jet-latitude index is known to exhibit a trimodal distribution, reflecting preferred southern, central, and northern jet positions that are commonly interpreted as distinct jet regimes. These regimes are typically defined by clustering the daily jet latitude values \citep{woollings2010variability, hannachi2017weather, strommen2023topological}. The southern (low jet latitudes), central (intermediate jet latitudes), and northern (high jet latitudes) regimes correspond to the negative, neutral, and positive phases of the North Atlantic Oscillation (NAO), respectively \citep{woollings2010variability, strommen2020jet}.

%====================================================
\section{Methods}\label{sec:method}
%====================================================

\subsection{ Persistence Homology}
\label{subsec:PH}
Persistent homology (PH) is one of the most widely used methods in Topological Data Analysis. It summarizes topological features—such as connected components, loops, and voids—across multiple spatial scales in a computationally accessible way.

The first step in PH is to construct a filtration, i.e., a nested family of "continuous shapes" built on the data. This is done using simple building blocks such as points, edges, triangles, and their higher-dimensional counterparts, called \emph{simplices}.

\begin{figure}[htbp]
    \centering
    \includegraphics[width=0.6\linewidth]{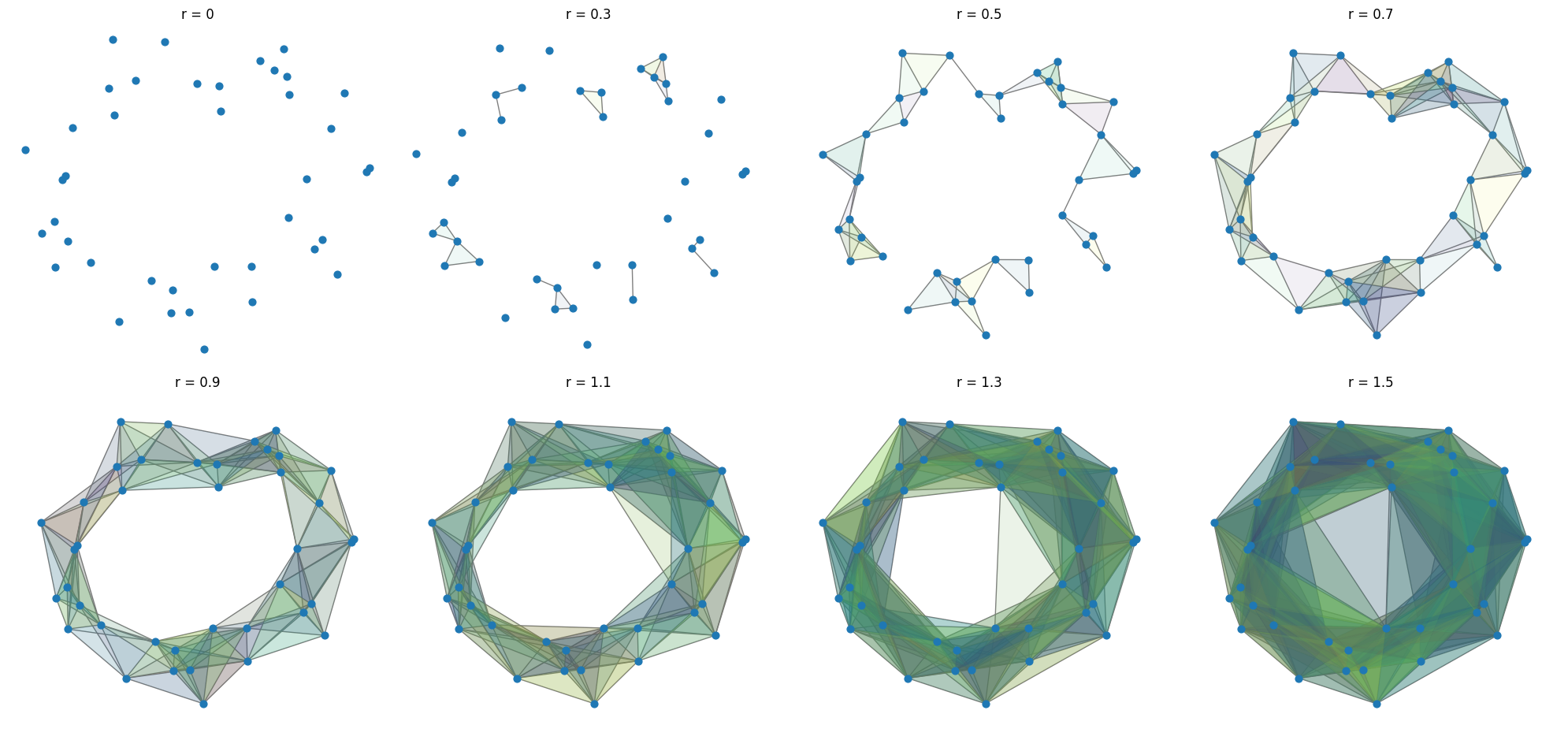}
   \caption{ Vietoris--Rips filtration of a noisy circle. For small radius values, edges connect only nearby points, forming small clusters. As the radius grows, these clusters merge and a loop appears. The loop later fills in as higher-dimensional simplices are added.}

    \label{fig:filtration}
\end{figure}

For example, nearby points are connected by edges whenever their distance is below a chosen threshold. Triangles and higher-dimensional simplices are added whenever all edges between their vertices are present. Varying the threshold therefore produces a sequence of nested shapes, called the Vietoris–Rips filtration, that captures the topology of the data across multiple spatial scales. Fig.~\ref{fig:filtration} illustrates this process for a noisy sample from a circle.

\begin{figure}[htbp]
    \centering
    \includegraphics[width=0.37\linewidth]{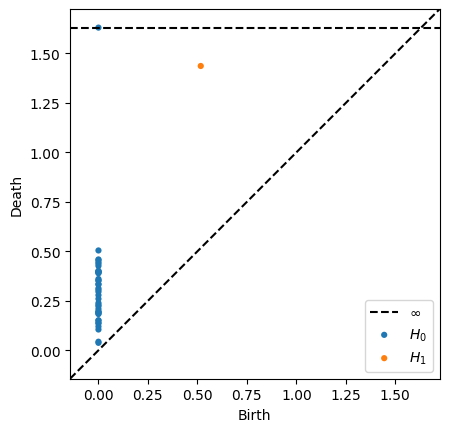}
    \caption{Persistence diagram showing two main topological features from the filtration in Fig.~\ref{fig:filtration}: the blue dot on the infinity line reflects the fact that all points eventually merge into a single connected component, while the orange dot corresponds to the lifespan of the main loop before the circle becomes filled in.
 }

    \label{fig:filtration_pd}
\end{figure}

PH records when each topological feature appears (its \emph{birth}) and when it disappears (its \emph{death}). Each feature is represented as a point $(\text{birth}, \text{death})$ in the plane, encoding its lifespan. The full collection of these points is called the \emph{persistence diagram} of the data. Because a feature must be born before it dies, all points lie on or above the diagonal. Points far from the diagonal correspond to features that persist across many scales and are often interpreted as meaningful structure, whereas points near the diagonal are typically less significant. Fig.~\ref{fig:filtration_pd} shows the persistence diagram for the filtration shown in Fig.~\ref{fig:filtration}. In this example, all points in the filtration eventually merge into a single connected component (which never disappears), while the prominent loop appears near radius~0.7 and fills in around radius~1.5. These features are visible in the corresponding persistence diagram: the long-lived connected component is shown in blue on the dashed line indicating infinity (a convention used for features that persist throughout the entire filtration), and the main loop is represented as $(0.7, 1.5)$ in orange.  

\noindent
For further background on persistent homology, we refer the reader to \citep{ Chazal2021TDA, edelsbrunner2002topological, ghrist2008barcodes, edelsbrunner2010computational, otter2017roadmap}.

%====================================================
\subsection{Centrality-Based Bifiltration}
%====================================================

In S23, a density–radius bifiltration was introduced to extract topological features across varying density thresholds. To better capture fine-scale features—such as the southern jet regime in JetLat and the thin loop structures in the CdV system—we replace density with a measure of centrality and introduce a centrality–radius bifiltration. For a dataset with $n$ points, we consider a family of continuous functions $\{C_k \mid k = 1, 2, \cdots, n\}$, each of which quantifies how central a point is relative to its local neighborhood at scale $k$.
 
These functions are derived from the well-known \emph{distance-to-measure} functions \citet{chazal2011geometric}. Given \(k\), the \(k\)-distance-to-measure ($k$-dtm) function \(d_k\) computes the average distance from a point to its $k$ nearest neighbors.

Formally, let $X = \{x_1, \ldots, x_n\} \subset \mathbb{R}^d$ be our dataset. The $k$-dtm is defined by

\[d_{k}(x) = \left( \frac{1}{k} \sum_{i=1}^{k} \| x - x_{(i)} \|^2 \right)^{1/2},
\]
where $x_{(1)}, \ldots, x_{(k)}$ are the $k$ nearest neighbors of $x$ under the Euclidean norm.

The motivation is straightforward: points in dense regions have many nearby neighbors and therefore take lower dtm values, while points in sparse or weakly connected regions have larger dtm values. 

The parameter \(k\) controls the scale of structure that is captured: small \(k\) values are more sensitive to fine structure and therefore to noise. On the other hand, large \(k\) values are less sensitive to small variations and emphasize more global geometric structure. Since
\[
d_1 \le d_2 \le \cdots \le d_n,
\]
increasing \(k\) can gradually control the level of sensitivity from very fine local geometry to the global organization of the data.

The \emph{local centrality function} $C_k$ used in this study is the normalized and inverted form of the $k$-dtm:
\[
C_k(x) = 1 - \frac{d_k(x) - \min_j d_k(x_j)}{\max_j d_k(x_j) - \min_j d_k(x_j)}.
\]

\begin{figure}[htbp]
    \centering
    \includegraphics[width=0.5\linewidth]{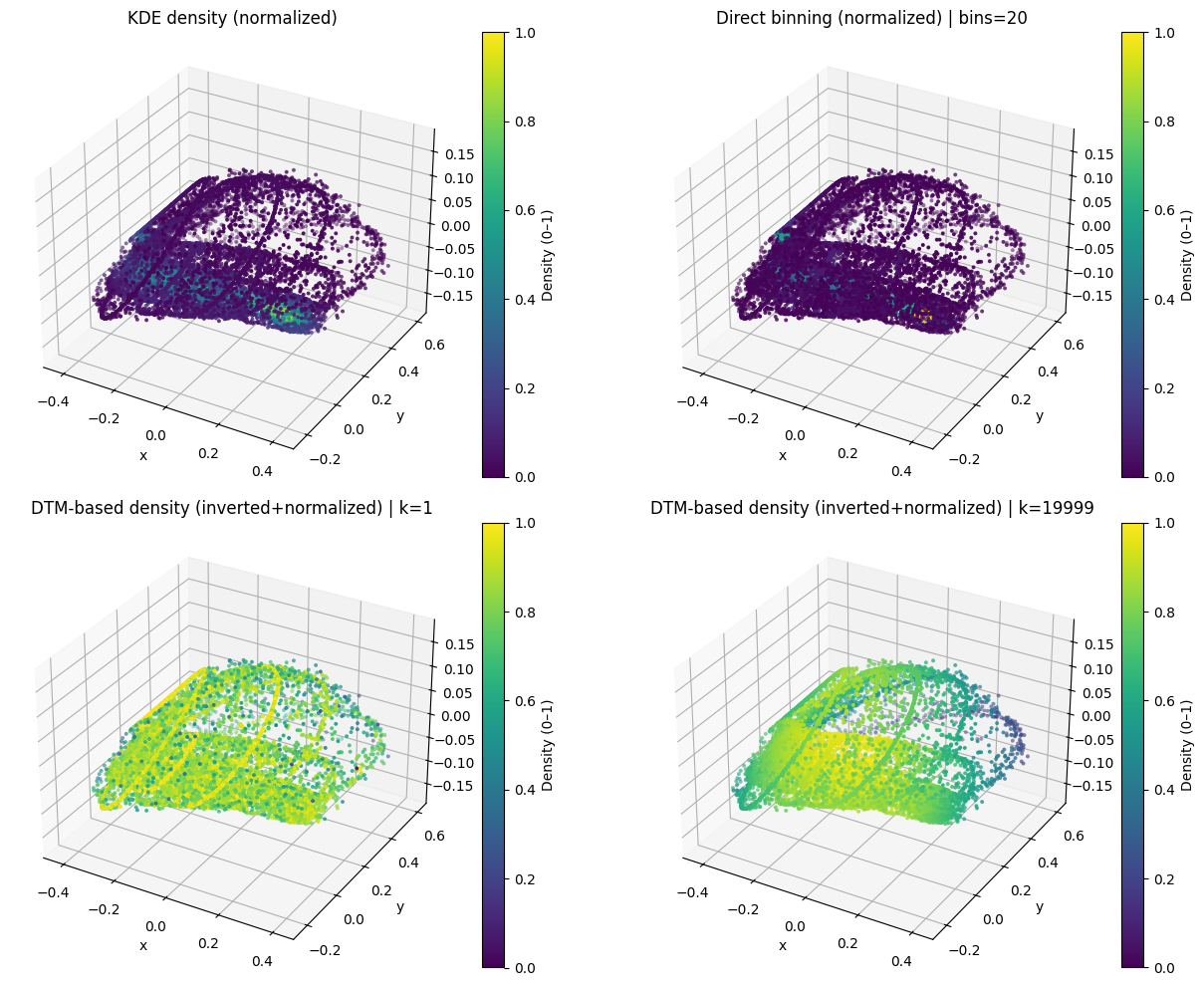}
    \caption{
CdV system with 20{,}000 data points colored by density (Gaussian KDE and direct binning) and centrality values.  Top-left: Gaussian KDE highlights the dense central region while assigning low values to the thin loop structures. Top-right: direct binning with 20 bins per coordinate, highlighting mainly the central region. Bottom-left: $C_1$ and bottom-right: $C_{19999}$ assign high values to both the central region and loop structures, with $C_1$ placing stronger emphasis on the loops.}
\label{fig:CDV_comparison}
\end{figure}

In practice, the choice of $k$ is empirical and depends on the dataset. A common starting point is to select $k$ as a small fraction of the dataset size, typically within the range $0.01\%$–$1\%$. For clean and well-sampled datasets, such as the toy models considered here, values of $k$ in this range perform well. As $k$ increases further, the centrality function becomes increasingly smooth and exhibits a global averaging behavior that qualitatively resembles KDE. Fig.~\ref{fig:CDV_comparison} compares KDE, direct binning, and the centrality functions ($C_1$ and $C_n$) on the CdV dataset, with each panel showing the same point cloud colored by the respective function values.

%==========================dtm%================
\subsection{Computational Methodology}\label{subsec:computation}
%====================================================

Our computational approach follows S23 with some extensions. It includes five main steps, summarized below.

\begin{enumerate}
    \item \textbf{Normalization.}  
    Each coordinate of the data set is normalized to have unit variance. The resulting data set is denoted by $X$.

    \item \textbf{Evaluation of Density or Centrality.}  
    We choose either Gaussian KDE, direct binning, or the centrality function \(C_k\) to assign a scalar value to each data point. Local densities are estimated using \texttt{scipy.stats.gaussian\_kde} with Scott’s rule for bandwidth selection. The distance-to-measure function \(d_k\) is computed using \(k\)-nearest neighbors (scikit-learn) via root-mean-square distances, and the centrality function \(C_k\) is obtained by min--max normalization followed by inversion.

    \item \textbf{Density or Centrality Thresholding.}  

    A filtration is constructed by selecting a sequence of percentiles $p_1 < p_2 < \dots < p_l$. For each percentile $p_i$, we form a subsample $X_{p_i}$ containing the top $p_i\%$ of points—those with the highest KDE values or, when using centrality, the highest centrality values.
     
    \item \textbf{Persistent Homology.}  
    Persistent homology is computed on each \(X_{p_i}\) using the \texttt{Gudhi} library. From each persistence diagram, we extract the longest-lived connected components and loops, along with their birth and death times. Representative cycles for the longest-lived loops are computed using \texttt{Persloop}. For clarity, only the five longest-lived features are retained; in our datasets, this does not remove any significant topological structure \citep{straus2007circulation, strommen2023topological}. The main computational settings follow those of S23.
    
    \item \textbf{Bifiltration.}  
    Steps (2)–(4) are repeated for \(P = 10\%, 20\%, \dots, 100\%\), producing a family of subsamples filtered by either density or centrality. These thresholds are sufficient for the datasets considered here, although finer discretizations may be useful for more complex dynamics.
\end{enumerate}

%==========================%Significance Testing%================
\subsection{Significance Testing and Reliability of Topological Features} \label{subsec:sig_test}
%====================================================

Not all features in a persistence diagram reflect meaningful structure; some arise from noise, sampling effects, or parameter choices. To identify reliable features, we adopted the significance tests of S23 and extended them to our centrality–based framework.

\medskip
\noindent
\textbf{(1) Sensitivity tests for persistent homology.}
We repeated the analysis using different filtering thresholds and homology settings. The main structure of the persistence diagrams and the dominant connected components remained stable, indicating robustness. Loop features were more variable in their persistence and are therefore interpreted as qualitative features of the attractor’s geometry rather than precise invariants.

\medskip
\noindent
\textbf{(2) Gaussian reference test.} 

To estimate the level of topological noise expected from an unstructured distribution, we generated $10{,}000$ points from a three-dimensional standard Gaussian distribution. We applied both the KDE–radius and centrality–radius bifiltrations, recording the \emph{maximum} lifespan of connected components and loops for each sample. The experiment was repeated ten times. The largest observed lifespan ($\leq 0.5$) was used as a noise threshold, and features exceeding this value were classified as non-trivial.

%==========================%Results%==========================
\section{Results}\label{sec:results}
%====================================================
The structure of each dataset is summarized using a \textit{bifiltration diagram}, which tracks how topological features persist as the density or centrality threshold changes. The horizontal axis shows the threshold value, and the vertical axis shows the lifespan of each feature. At each threshold, the five most persistent features are displayed.

The dashed line indicates the maximum lifespan obtained from the Gaussian reference test and serves as a baseline for topological noise, while the \texttt{min-pers} line marks the minimum lifespan required for a feature to be considered significance. Connected components with fewer than four points are treated as noise (red “×”), those with four to ten points as weak features (yellow diamonds), and those with more than ten points as robust features (red circles). Loop features are shown as blue triangles (see Fig.~\ref{fig:L63_result}(a) for an example).

%=========================%L63%==========================%
\subsection{Results for the Lorenz--63 System}\label{subsec:l63result}
%====================================================
Recall that the Lorenz--63 attractor resembles a “butterfly” shape, with each wing commonly interpreted as a distinct regime. In the bifiltration framework, these regimes are associated with two persistent loops in the bifiltration diagram of the attractor.

\begin{figure}[H]
    \centering
    % -------- Row 1 --------
    \begin{subfigure}{0.23\textwidth}
        \includegraphics[width=\linewidth]{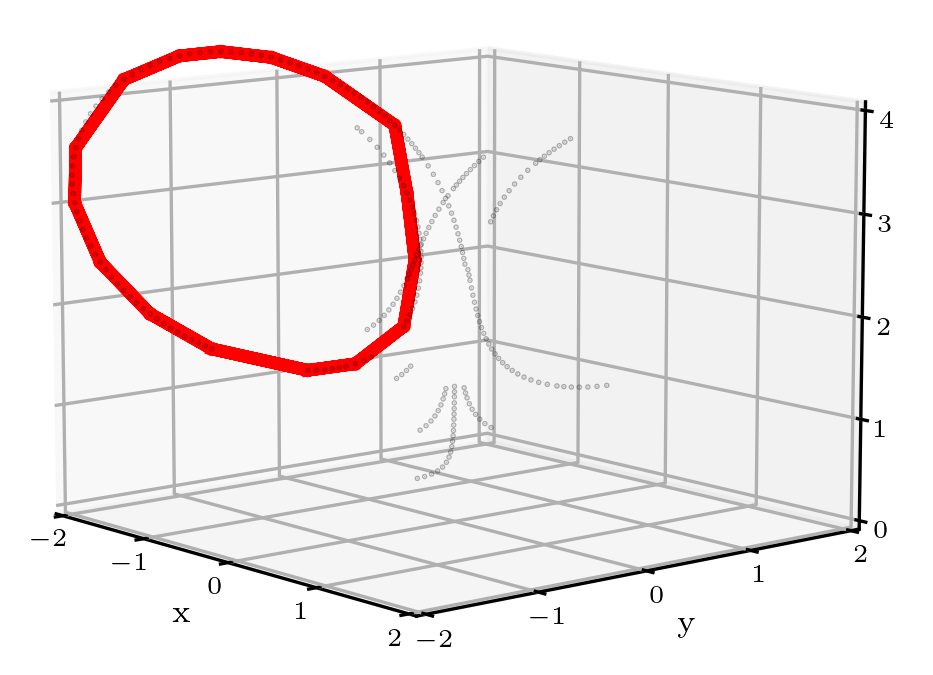}
        \caption{$C_1$: 10\%}
        \label{fig:L63_loops_C1}
    \end{subfigure}
    \begin{subfigure}{0.23\textwidth}
        \centering
        \includegraphics[width=\linewidth]{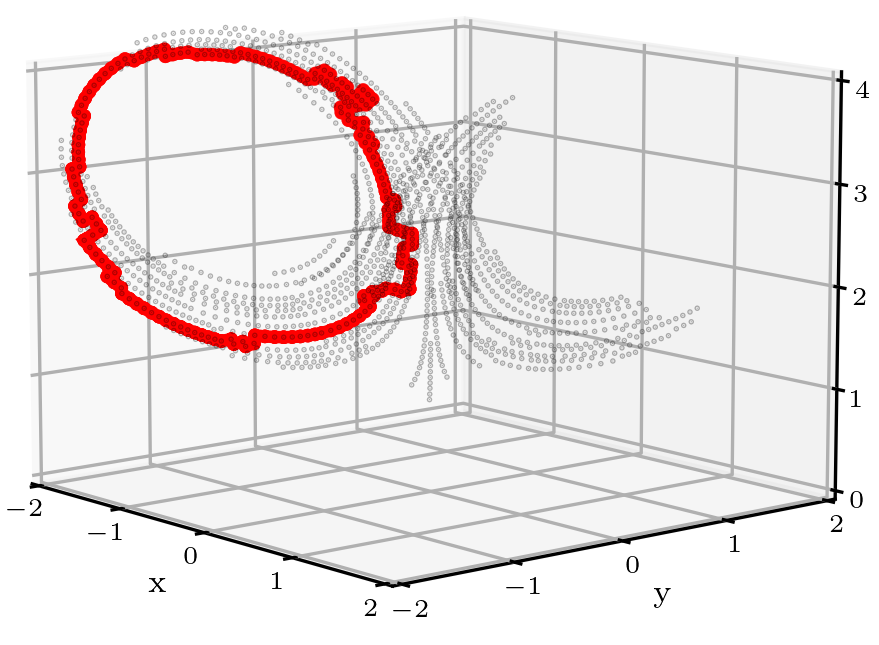}
        \caption{$C_{100}$: 60\%}
        \label{fig:L63_loops_C100}
    \end{subfigure}
    % -------- Row 2 --------
    \begin{subfigure}{0.23\textwidth}
        \includegraphics[width=\linewidth]{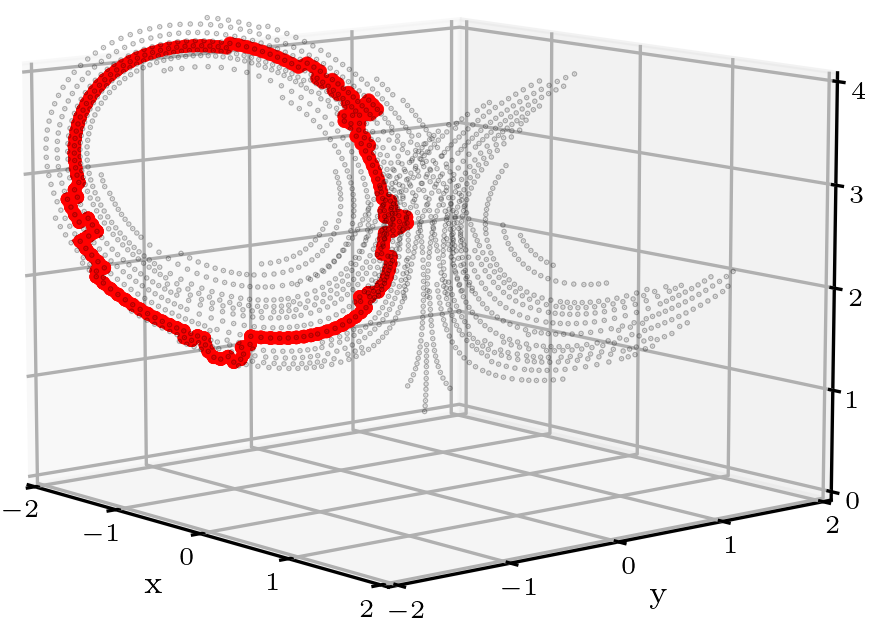}
        \caption{$C_n$: 70\%}
        \label{fig:L63_loops_C200}
    \end{subfigure}
    \begin{subfigure}{0.23\textwidth}
        \includegraphics[width=\linewidth]{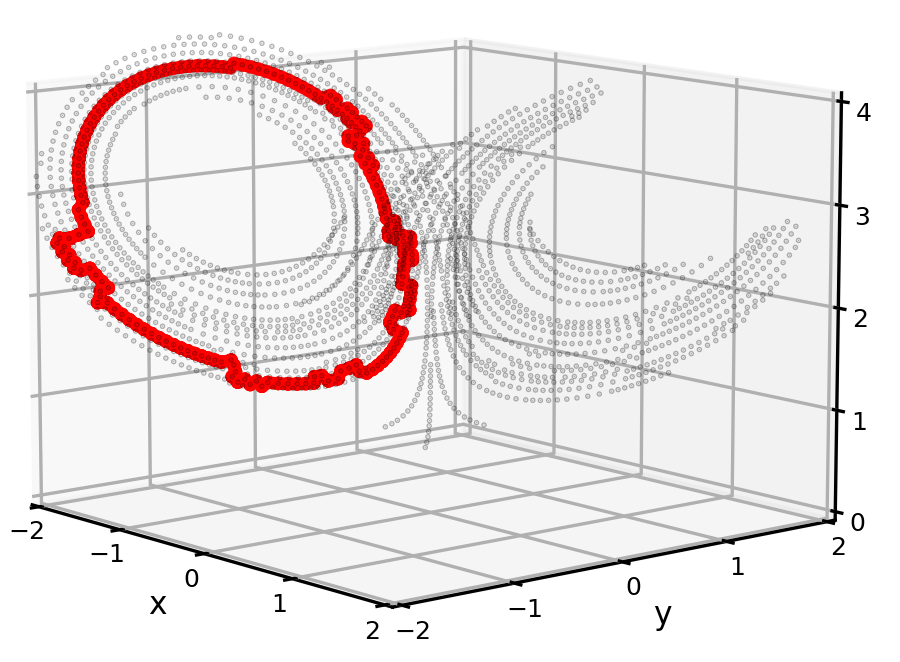}
        \caption{KDE: 80\%}
        \label{fig:L63_loops_kde}
    \end{subfigure}

    \caption{Lorenz--63: the first persistent loop (left lobe) identified by both methods. Panel (c) shows \(C_n\), where \(n\) denotes the full neighborhood size.}
    \label{fig:L63_loops}
\end{figure}

\begin{figure}[H]
    \centering

    % ===================== Row 1 =====================
    \begin{subfigure}[t]{0.3\textwidth}
        \centering
        \includegraphics[width=\linewidth]{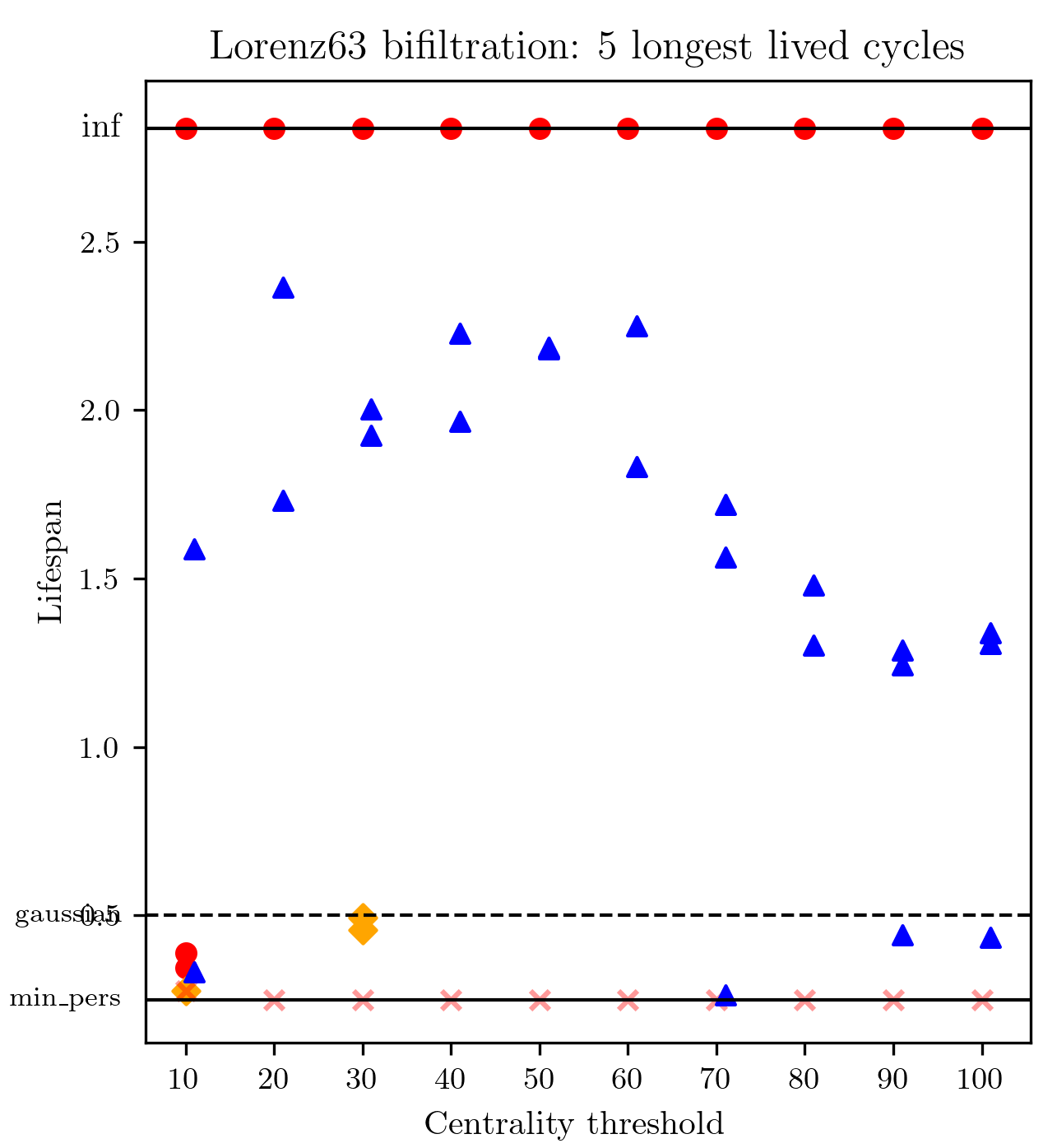}
        \caption{$C_1$-radius bifiltration diagram}
        \label{fig:L63_C1}
    \end{subfigure}
    \hspace{0.4cm}
    \begin{subfigure}[t]{0.3\textwidth}
        \centering
        \includegraphics[width=\linewidth]{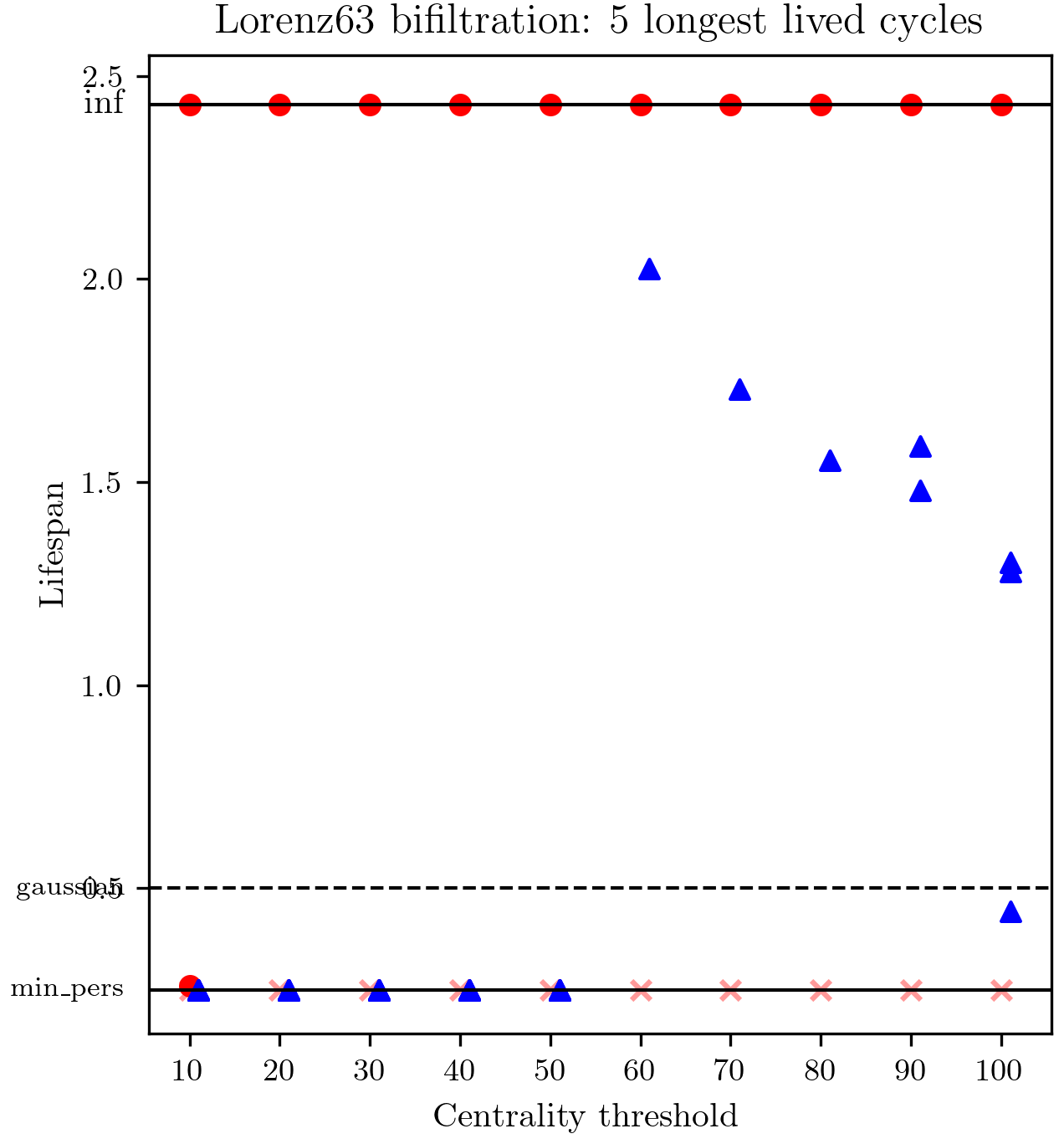}
        \caption{$C_{100}$-radius bifiltration diagram}
        \label{fig:L63_C100}
    \end{subfigure}

    \vspace{0.5em}

    % ===================== Row 2 =====================
    \begin{subfigure}[t]{0.3\textwidth}
        \centering
        \includegraphics[width=\linewidth]{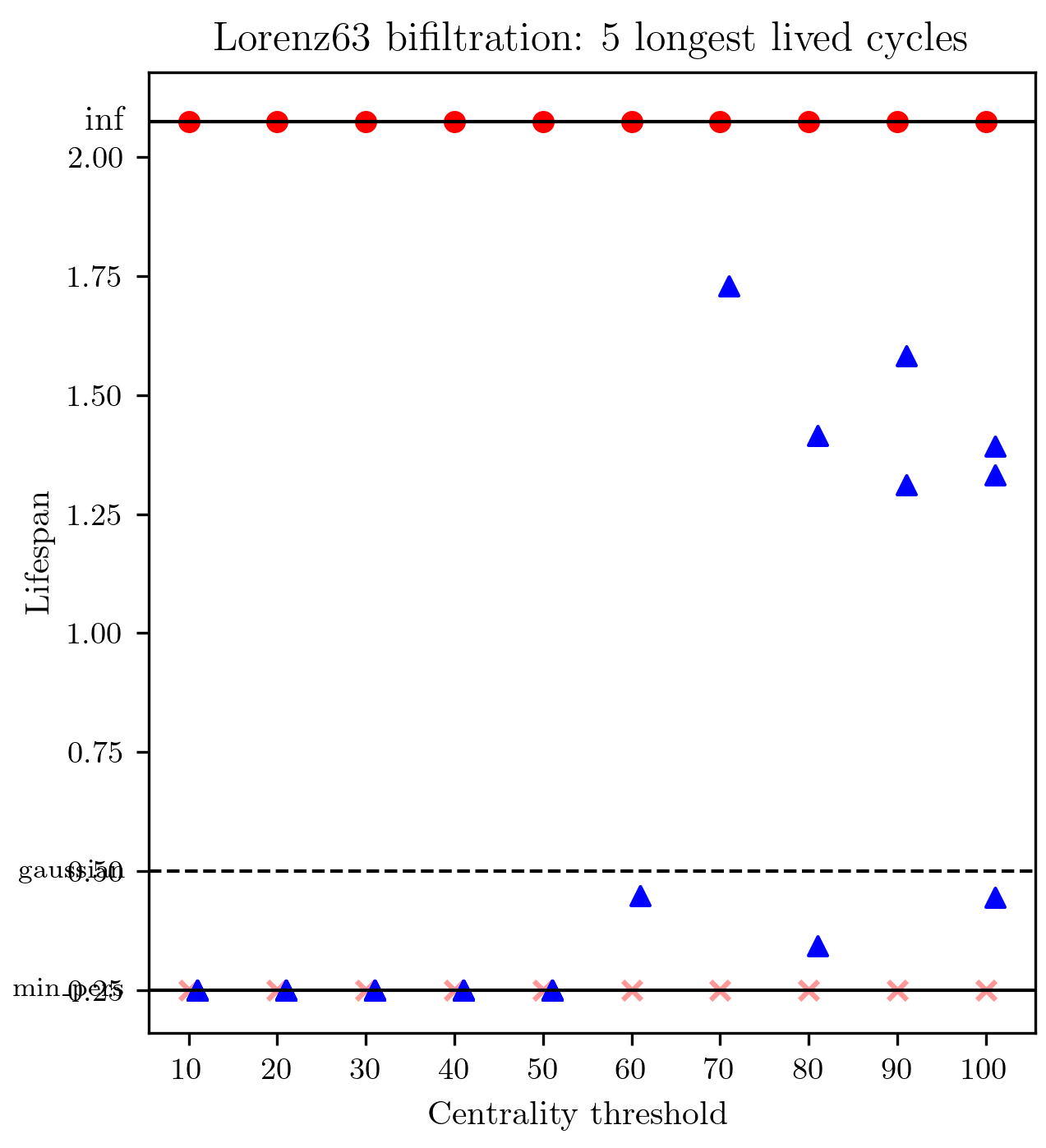}
        \caption{$C_n$-radius bifiltration diagram (\(n\) denotes the full neighborhood size)}
        \label{fig:L63_Cn}
    \end{subfigure}
    \hspace{0.4cm}
    \begin{subfigure}[t]{0.3\textwidth}
        \centering
        \includegraphics[width=\linewidth]{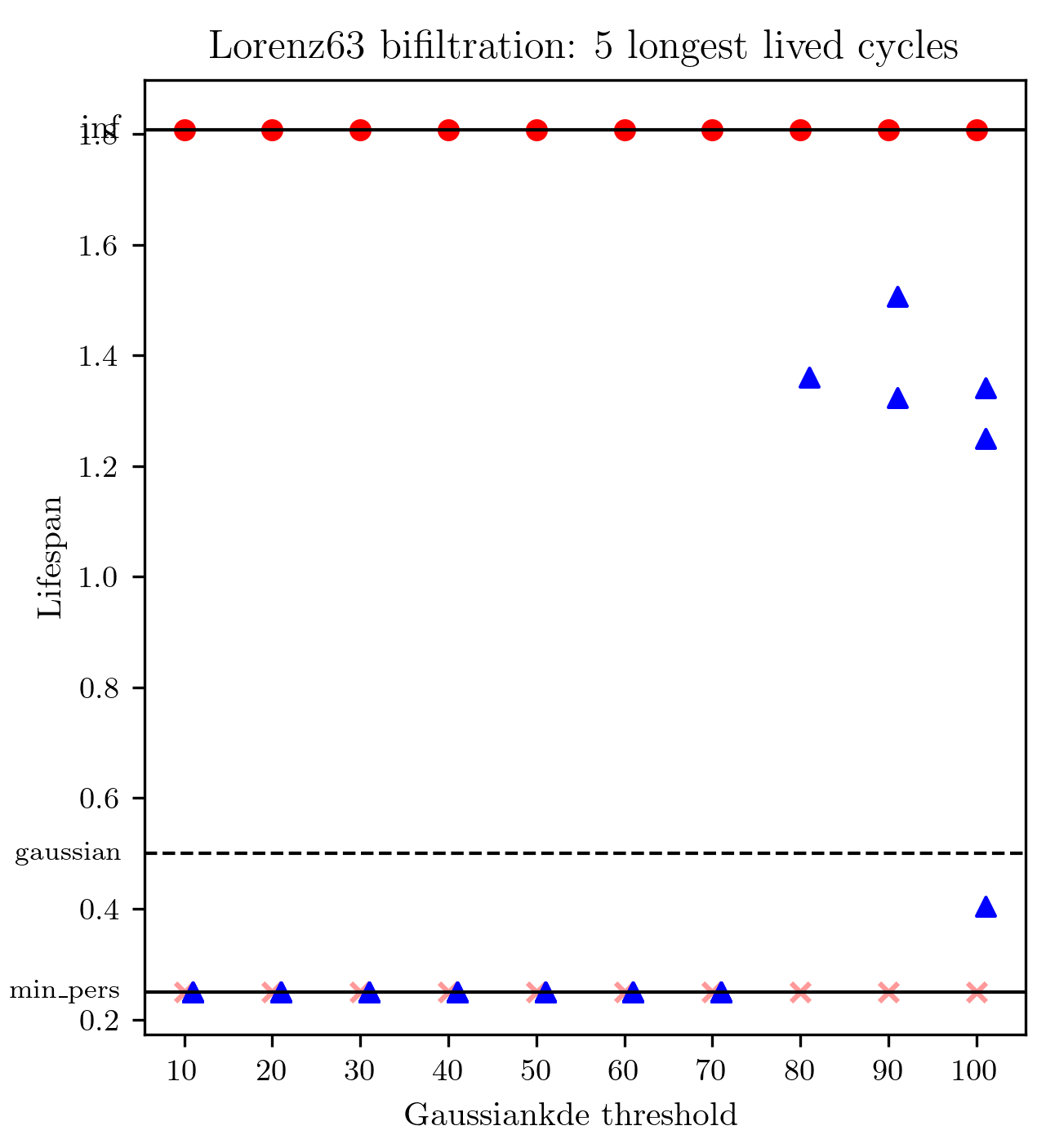}
        \caption{Gaussian KDE-radius bifiltration diagram}
        \label{fig:L63_KDE}
    \end{subfigure}

    \caption{
    Lorenz--63: comparison of centrality-based and KDE-based bifiltrations. Features below the Gaussian noise level are treated as noise.
    }

    \label{fig:L63_result}
\end{figure}

 Both the KDE-based and centrality-based methods (across all values of \(k\)) recover these loops. However, for smaller values of \(k\) (\(k \leq 100\)), the centrality-based method detects the corresponding loops at lower threshold levels and across a wider range of thresholds. Fig.~\ref{fig:L63_loops} shows the persistent loop corresponding to the left lobe, together with the threshold at which it is first detected by each method. For larger values of \(k\), the behavior of the centrality-based method becomes increasingly similar to that of the KDE-based approach. This pattern is illustrated in Fig.~\ref{fig:L63_result}.

%==========================%CDV%==========================%
\subsection*{Results for the Charney--DeVore System}\label{subsec:cdv}
%====================================================

Recall that the CdV system has two regimes corresponding to blocked and zonal states. In phase space, the blocked regime forms a dense, connected region, while the zonal regime appears as thin, low-density loops spread over a broader region.

\begin{figure}[H]
    \centering

    % ===================== Row 1 =====================
    \begin{subfigure}[t]{0.3\textwidth}
        \centering
        \includegraphics[width=\linewidth]{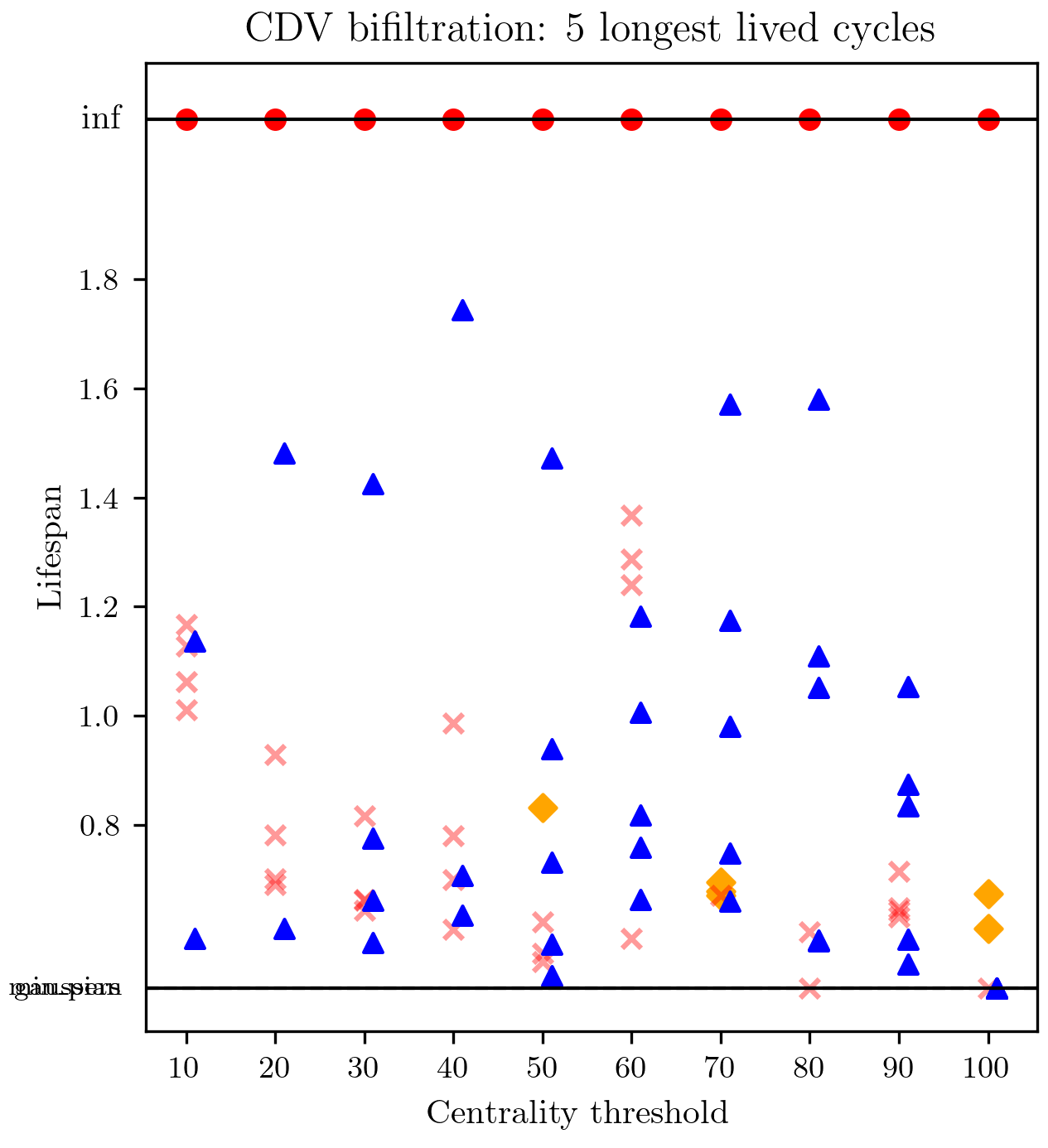}
        \caption{$C_1$-radius bifiltration diagram.}
        \label{fig:result_CDV_C1}
    \end{subfigure}
    \hspace{0.5cm}
    \begin{subfigure}[t]{0.3\textwidth}
        \centering
        \includegraphics[width=\linewidth]{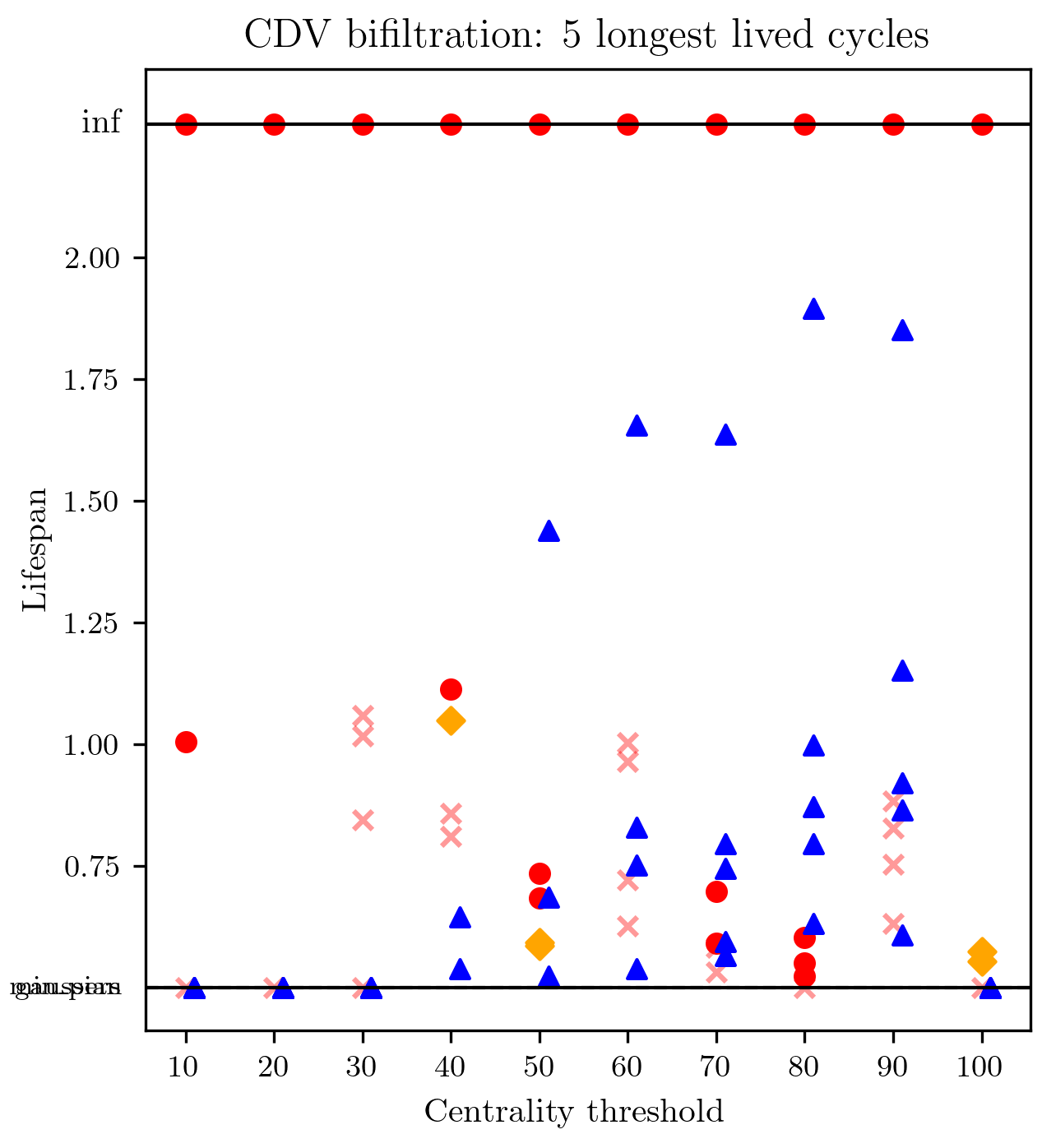}
        \caption{$C_{50}$-radius bifiltration diagram.}
        \label{fig:CDV_C50}
    \end{subfigure}

    \vspace{0.5em}

    % ===================== Row 2 =====================
    \begin{subfigure}[t]{0.3\textwidth}
        \centering
        \includegraphics[width=\linewidth]{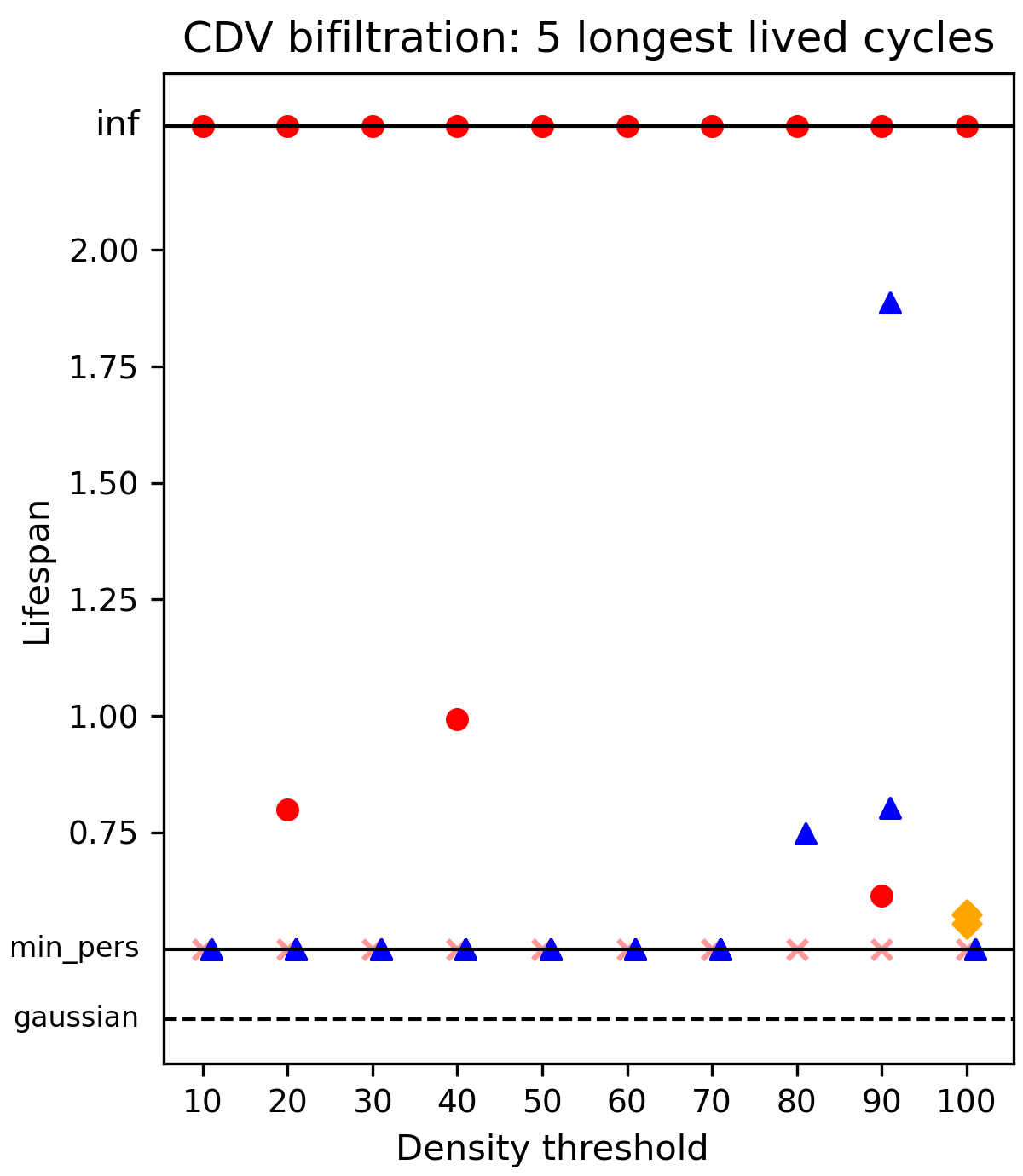}
        \caption{Gaussian KDE-radius bifiltration diagram.}
        \label{fig:CDV_result_KDE}
    \end{subfigure}
    \hspace{0.5cm}
    \begin{subfigure}[t]{0.3\textwidth}
        \centering
        \includegraphics[width=\linewidth]{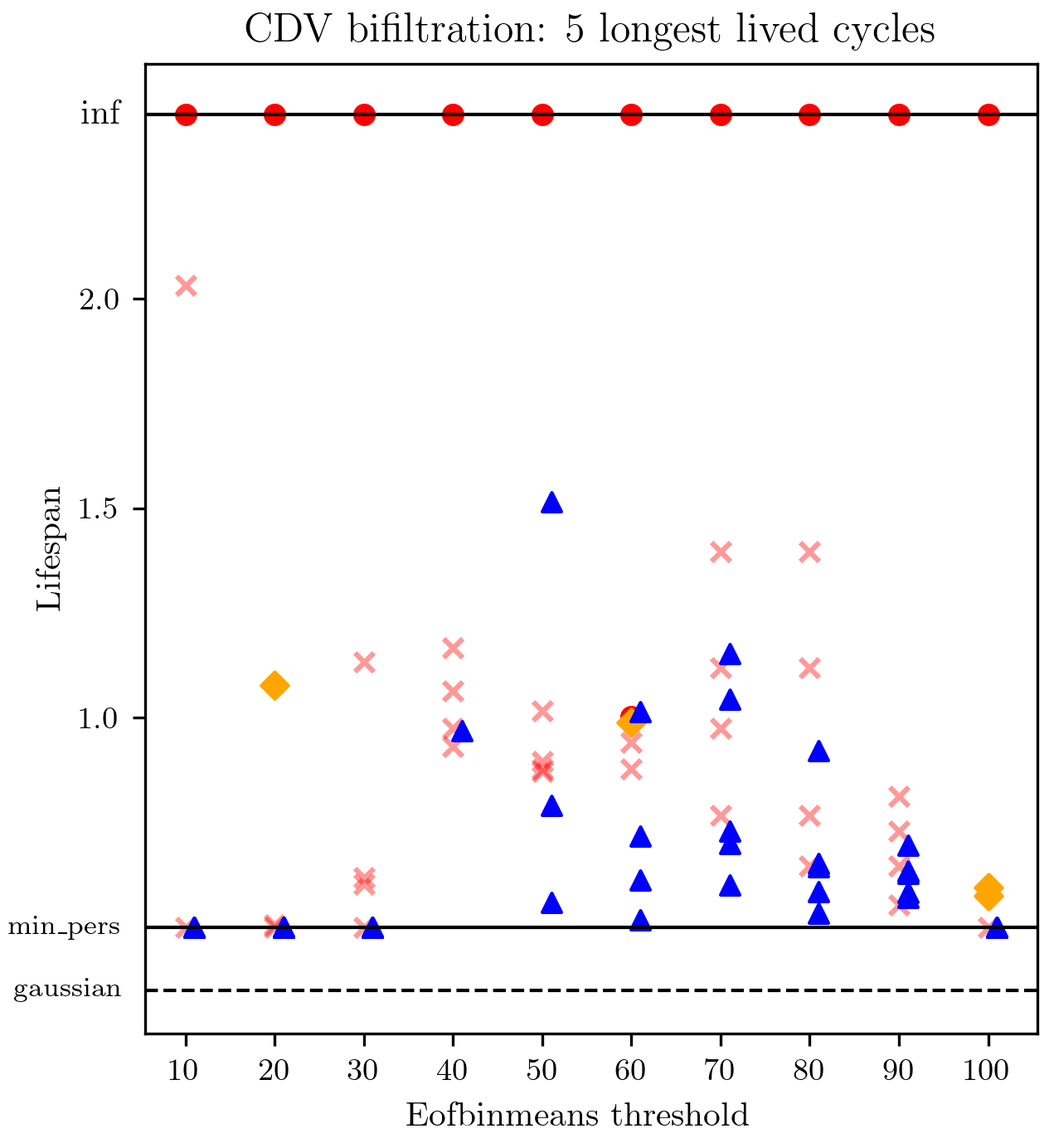}
        \caption{Binning-based density-radius bifiltration diagram.}
        \label{fig:CDV_Binning}
    \end{subfigure}

    \caption{
    CdV system: comparison of centrality-based, KDE-based, and binning-based bifiltrations. The KDE-based method suppresses thin loop structures associated with the zonal regime at low density thresholds. The centrality-based method detects these structures earlier in the filtration and produces more persistent features than the binning-based approach.
    }
    \label{fig:CDV_results}
\end{figure}

As shown in Fig.~\ref{fig:CDV_comparison}, Gaussian KDE tends to assign low density values to these thin loops, making them less likely to persist during the bifiltration process. This is reflected in the KDE-based bifiltration diagram (Fig.~\ref{fig:CDV_result_KDE}), where loop structures associated with the zonal regime are not reliably detected, particularly at low density thresholds. 

In S23, a direct-binning approach was introduced as an alternative to KDE. Compared with Gaussian KDE, it preserves more fine-scale loop structures, but some thinner loops are still weak or missed.

\begin{figure}[H]
  \centering
  % ===================== Row 1 =====================
  \begin{subfigure}[t]{0.35\textwidth}
    \centering
    \includegraphics[width=\linewidth]{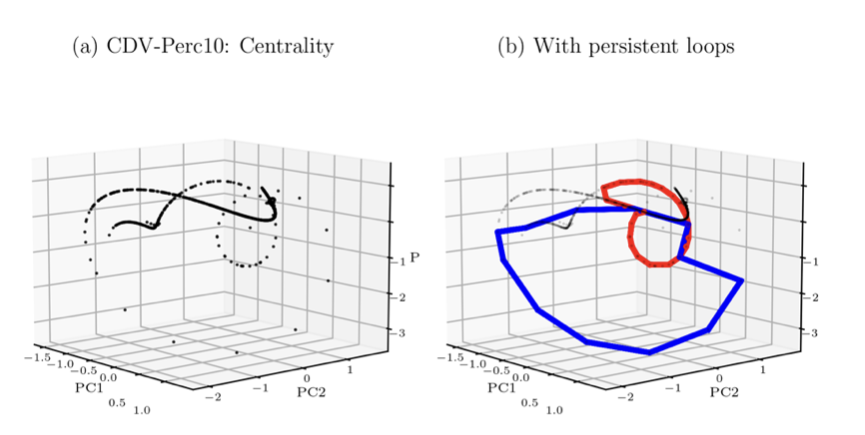}
    \caption{$C_1$: 10\%.}
    \label{fig:CDV_loops_C1}
  \end{subfigure}
  \hspace{0.1cm}
  \begin{subfigure}[t]{0.35\textwidth}
    \centering
    \includegraphics[width=\linewidth]{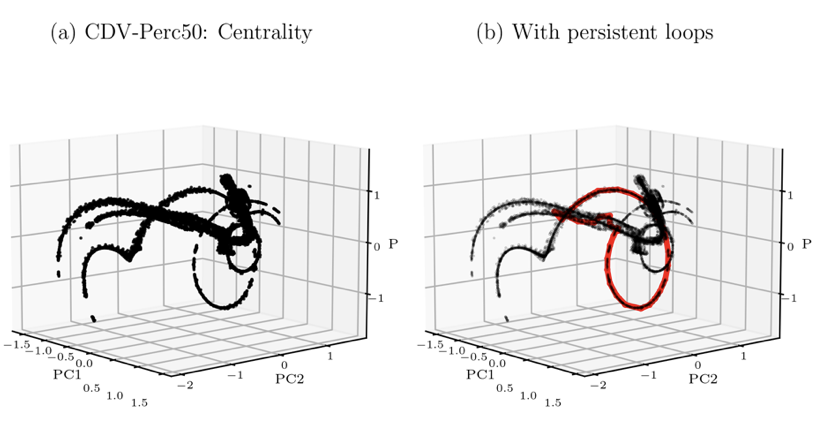}
    \caption{$C_{50}$: 50\%.}
    \label{fig:CDV_loops_C50}
  \end{subfigure}

  \vspace{0.5em}

  % ===================== Row 2 =====================
  \begin{subfigure}[t]{0.35\textwidth}
    \centering
    \includegraphics[width=\linewidth]{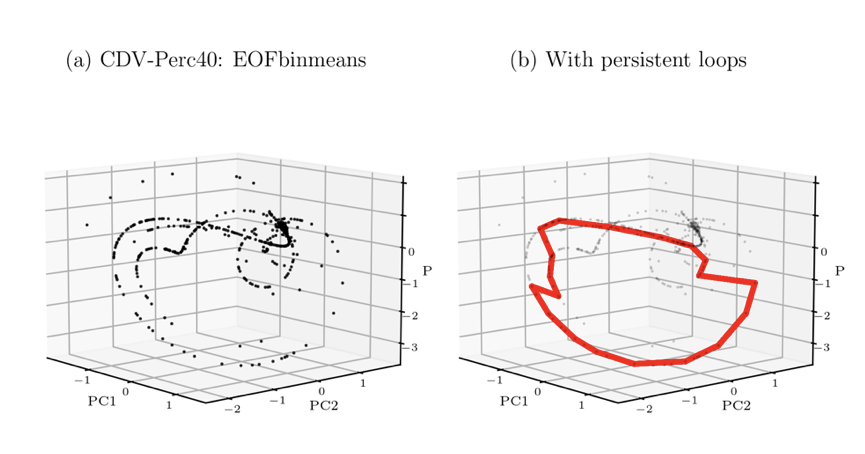}
    \caption{Binning-based density: 40\%.}
    \label{fig:CDV_loops_bins}
  \end{subfigure}
  \hspace{0.1cm}
  \begin{subfigure}[t]{0.35\textwidth}
    \centering
    \includegraphics[width=\linewidth]{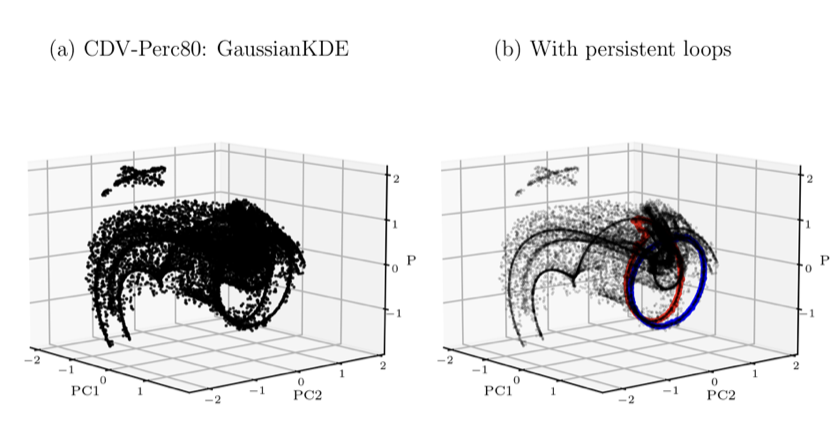}
    \caption{Gaussian KDE: 80\%.}
    \label{fig:CDV_loops_KDE}
  \end{subfigure}

\caption{
CdV dataset: first persistent loops detected by different bifiltration methods. The \(C_1\)-based method detects two loops at the 10\% threshold, \(C_{50}\) detects the first loop at 50\%, the binning-based method at 40\%, and Gaussian KDE at 80\%. Colors (assigned by \texttt{persloop}) indicate persistence, with red denoting the longest-lived loop.
}
\label{fig:CDV_loops}
\end{figure}

By contrast, the centrality-based method detects more loop structures and identifies them earlier in the filtration. The resulting features also exhibit longer persistence than those obtained from either the KDE-based or binning-based approaches. Fig.~\ref{fig:CDV_results} compares the corresponding bifiltration diagrams, while Fig.~\ref{fig:CDV_loops} shows the earliest persistent loops identified by each method, along with the threshold at which they are first detected.

%==========================%JetLat%==========================%
\subsection*{Results for the JetLat Dataset}\label{subsec:jetlat}
%====================================================

Recall that the JetLat dataset provides a three-dimensional representation of North Atlantic jet variability, consisting of the daily jet latitude index together with the first two principal components of 850 hPa zonal wind anomalies. 

\begin{figure}[H]
    \centering
    \begin{subfigure}[t]{0.3\linewidth}
        \centering
        \includegraphics[width=\linewidth]{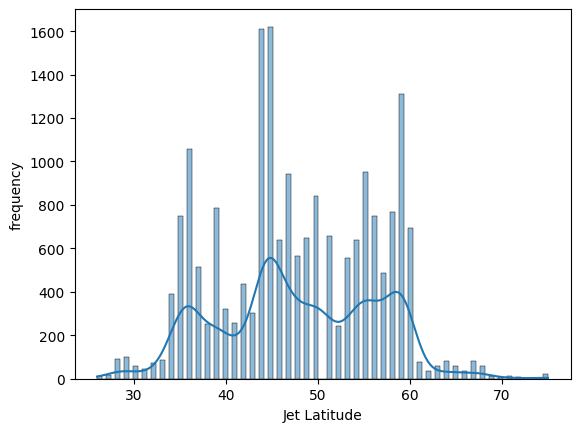}
        \caption{Distribution of jet latitude values in the JetLat dataset.}
        \label{fig:jet_distribution}
    \end{subfigure}
    \hspace{0.5 cm}
    \begin{subfigure}[t]{0.3\linewidth}
        \centering
        \includegraphics[width=\linewidth]{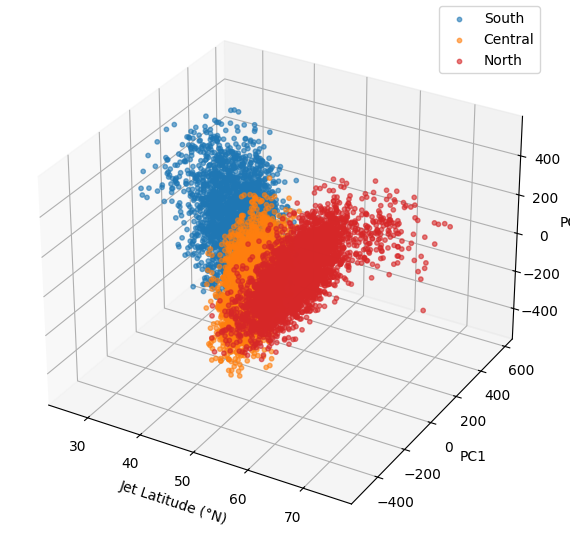}
        \caption{Organization of Southern, Central, and Northern jet regimes.}
    
        \label{fig:jet_regime}
    \end{subfigure}

    \caption{Jet latitude distribution and regime organization in the JetLat dataset. Regime assignment is based on local minima in the jet latitude distribution.}
    \label{fig:jetlat_regimes}
\end{figure}

The jet latitude index shows a distribution with three preferred states (southern, central, and northern jet positions), which are commonly interpreted as distinct circulation regimes. Fig.~\ref{fig:jetlat_regimes} shows the configuration of these regimes in the JetLat dataset. Regime boundaries are defined by local minima of the empirical jet-latitude distribution, which separate the preferred jet positions.

\begin{figure}[H]
    \centering
    % ===================== Row 1 =====================
    \begin{subfigure}[t]{0.3\textwidth}
        \centering
        \includegraphics[width=\linewidth]{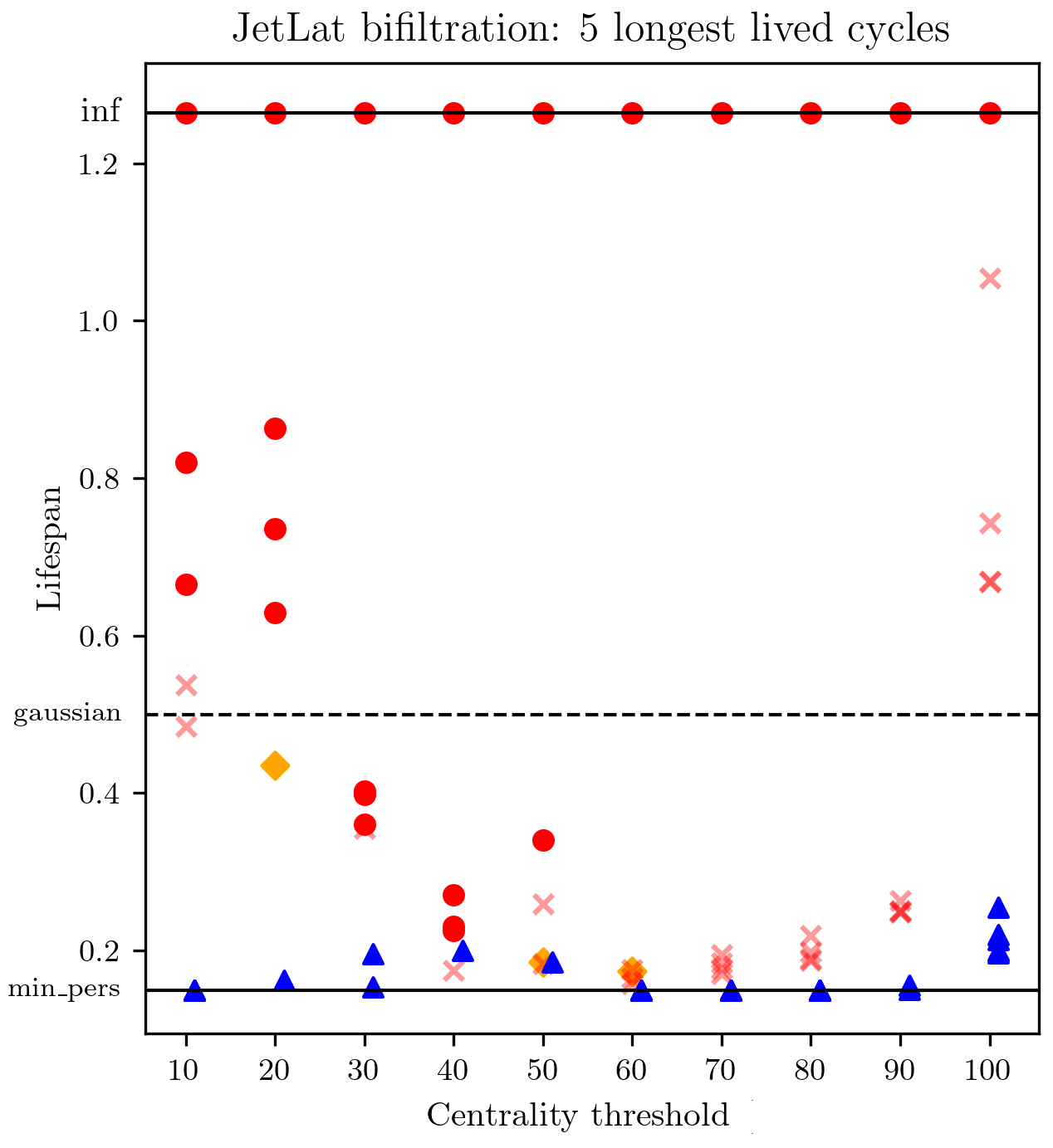}
        \caption{$C_{40}$-radius bifiltration diagram.}
        \label{fig:JetLat_C40}
    \end{subfigure}
    \hspace{0.5cm}
    \begin{subfigure}[t]{0.3\textwidth}
        \centering
        \includegraphics[width=\linewidth]{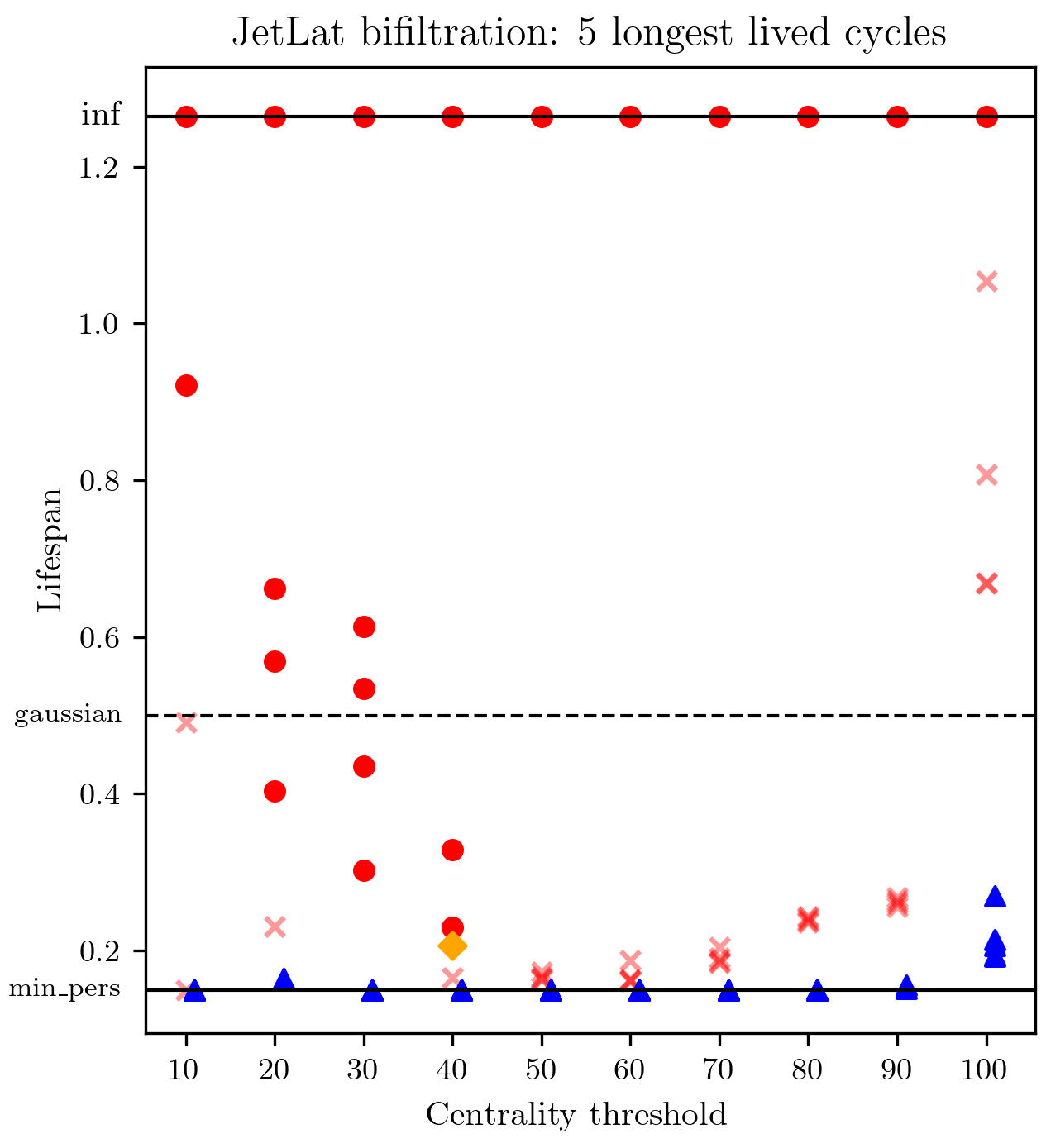}
        \caption{$C_{75}$-radius bifiltration diagram.}
        \label{fig:JetLat_C75}
    \end{subfigure}

    \vspace{0.5em}

    % ===================== Row 2 =====================
    \begin{subfigure}[t]{0.3\textwidth}
        \centering
        \includegraphics[width=\linewidth]{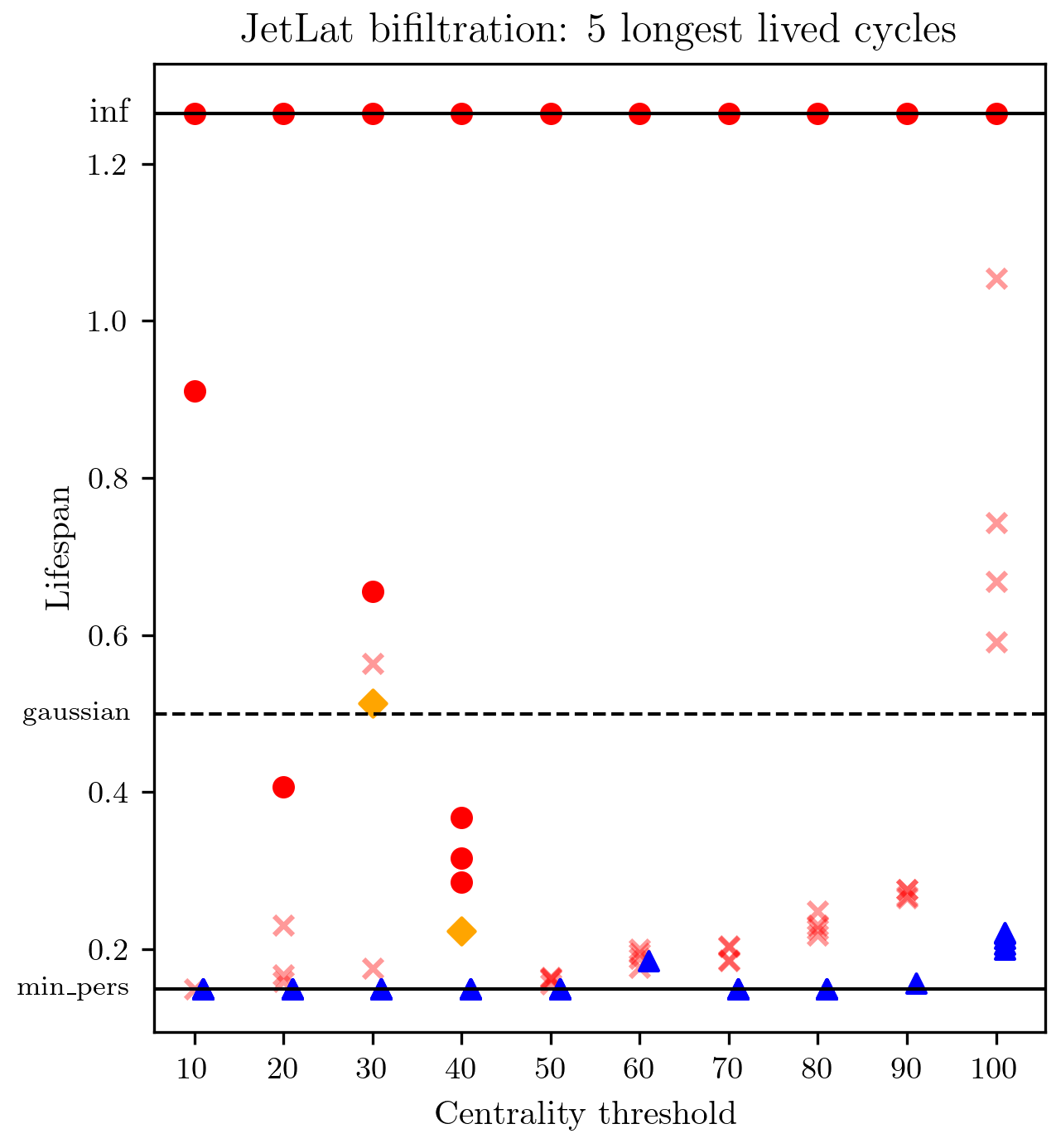}
        \caption{$C_{100}$-radius bifiltration diagram.}
        \label{fig:JetLat_C100}
    \end{subfigure}
    \hspace{0.5cm}
    \begin{subfigure}[t]{0.3\textwidth}
        \centering
        \includegraphics[width=\linewidth]{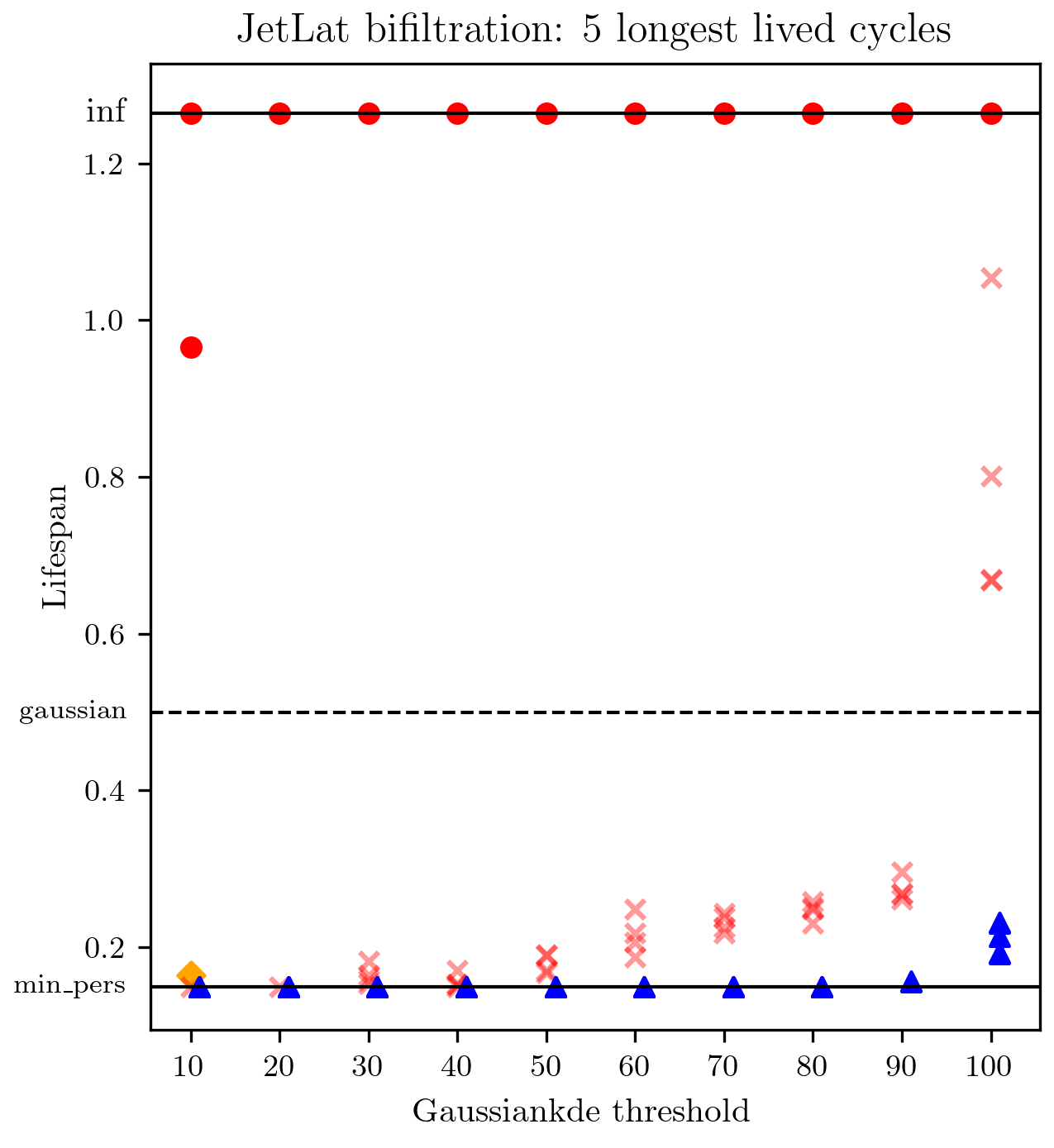}
        \caption{Gaussian KDE-radius bifiltration diagram.}
        \label{fig:JetLat_result_KDE}
    \end{subfigure}

   \caption{
    JetLat dataset: comparison of centrality-based and KDE-based bifiltrations. Increasing the locality scale of the centrality functions ($C_{40}$, $C_{75}$, $C_{100}$) suppresses fine-scale regime structure and yields behavior closer to the KDE-based bifiltration.
    }
    \label{fig:JetLat_results}
\end{figure}

The Gaussian KDE-based bifiltration fails to recover all three jet regimes. At the 10\% density threshold, it identifies only two connected components, corresponding to the central and northern jet regimes (Fig.~\ref{fig:JetLat_result_KDE}). The southern jet regime does not appear as a separate component at any threshold. As noted by S23, neither refinement of the density threshold nor increasing the number of principal components improves this outcome, suggesting that the two-regime pattern is a robust outcome of the Gaussian KDE-based approach.

On the other hand, the centrality-based method yields a richer structure. At low centrality thresholds ($p \leq 40\%$), it consistently identifies two connected components across all tested $k$. For intermediate neighborhood sizes ($30 \leq k \leq 100$), this number increases to three or four components at low $p$, while outside this range of $k$ it again collapses to two (Fig.~\ref{fig:JetLat_results}).

To understand the nature of these components, we assign each one to a regime based on the dominant jet-latitude proportion among its points. Some components are strongly associated with a single regime, while others mix regimes; we refer to these as mixed-regime components. To quantify this, we define regime purity as the fraction of points in a component belonging to its assigned dominant regime. Table~\ref{table:regime_assignment} summarizes results for $k=40$.

\begin{table}[H]
\centering
\begin{tabular}{c r l r l}
\hline
$p$ (\%) & component size & assigned regime & purity & regime counts $(S,C,N)$ \\
\hline
10 & 287  & North   & 0.951 & $(2,12,273)$ \\
10 & 676  & Central & 0.964 & $(10,652,14)$ \\
10 & 13   & South   & 0.769 & $(10,3,0)$ \\
\hline
20 & 1902 & Central & 0.640 & $(26,1217,659)$ \\
20 & 20   & South   & 0.700 & $(14,6,0)$ \\
20 & 12   & South   & 1.000 & $(12,0,0)$ \\
20 & 19   & South   & 1.000 & $(19,0,0)$ \\
\hline
\end{tabular}
\caption{Regime assignment for connected components identified by $C_{40}$-radius bifiltration.}
\label{table:regime_assignment}
\end{table}

The presence of these mixed and fragmented structures is consistent with previous analyses of North Atlantic jet variability. In particular, \citet{madonna2017link} identify a mixed jet cluster that contains distinct configurations, including split and strongly tilted jets, which are associated with European and Scandinavian blocking.

\begin{figure}
\centering

% ===================== Row 1: C40 =====================
\begin{subfigure}[t]{0.25\textwidth}
    \centering
    \includegraphics[width=\linewidth]{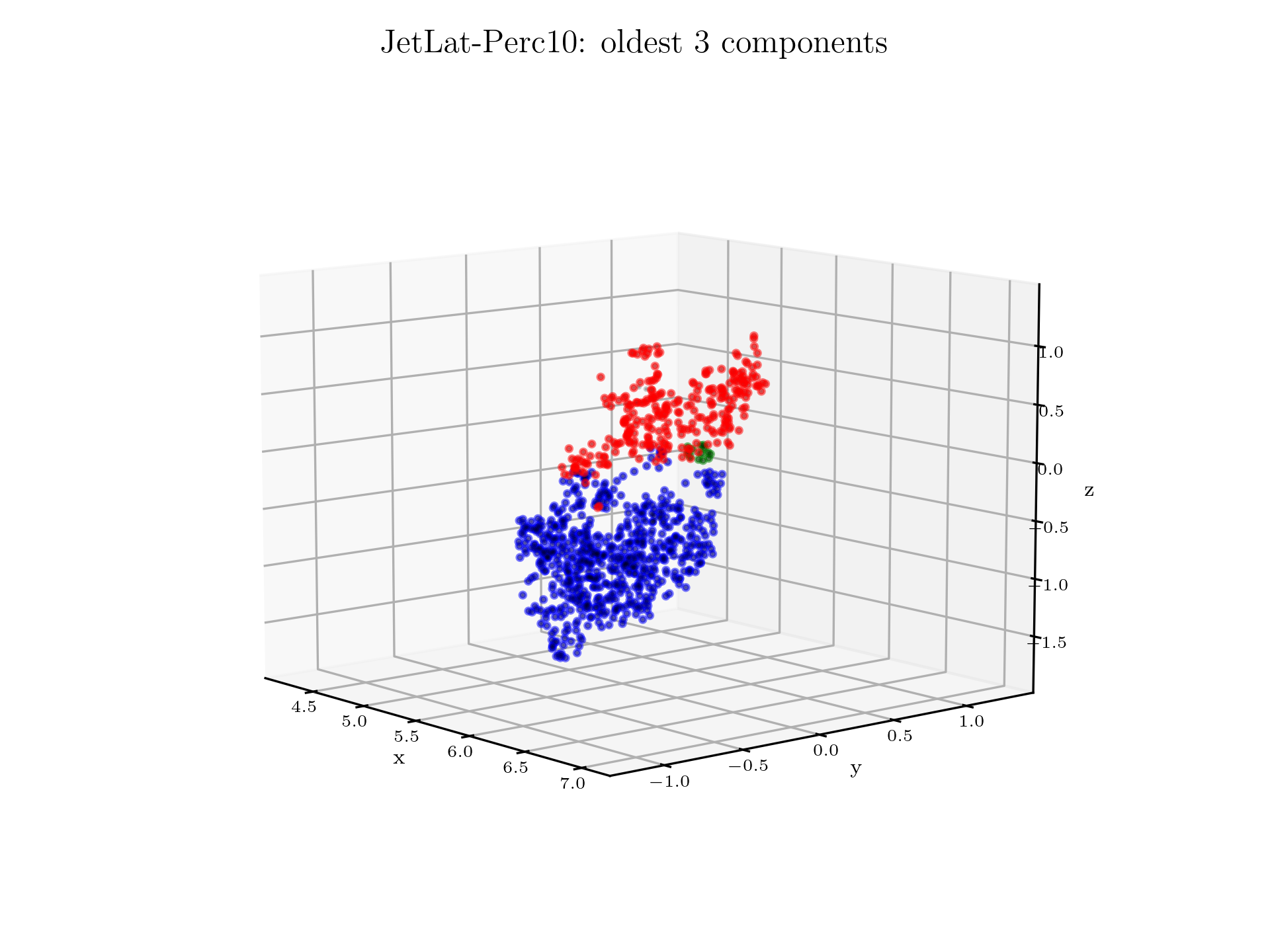}
    \caption{$C_{40}$ top 10\%}
    \label{fig:Jet_C40_10}
\end{subfigure}\hfill
\begin{subfigure}[t]{0.25\textwidth}
    \centering
    \includegraphics[width=\linewidth]{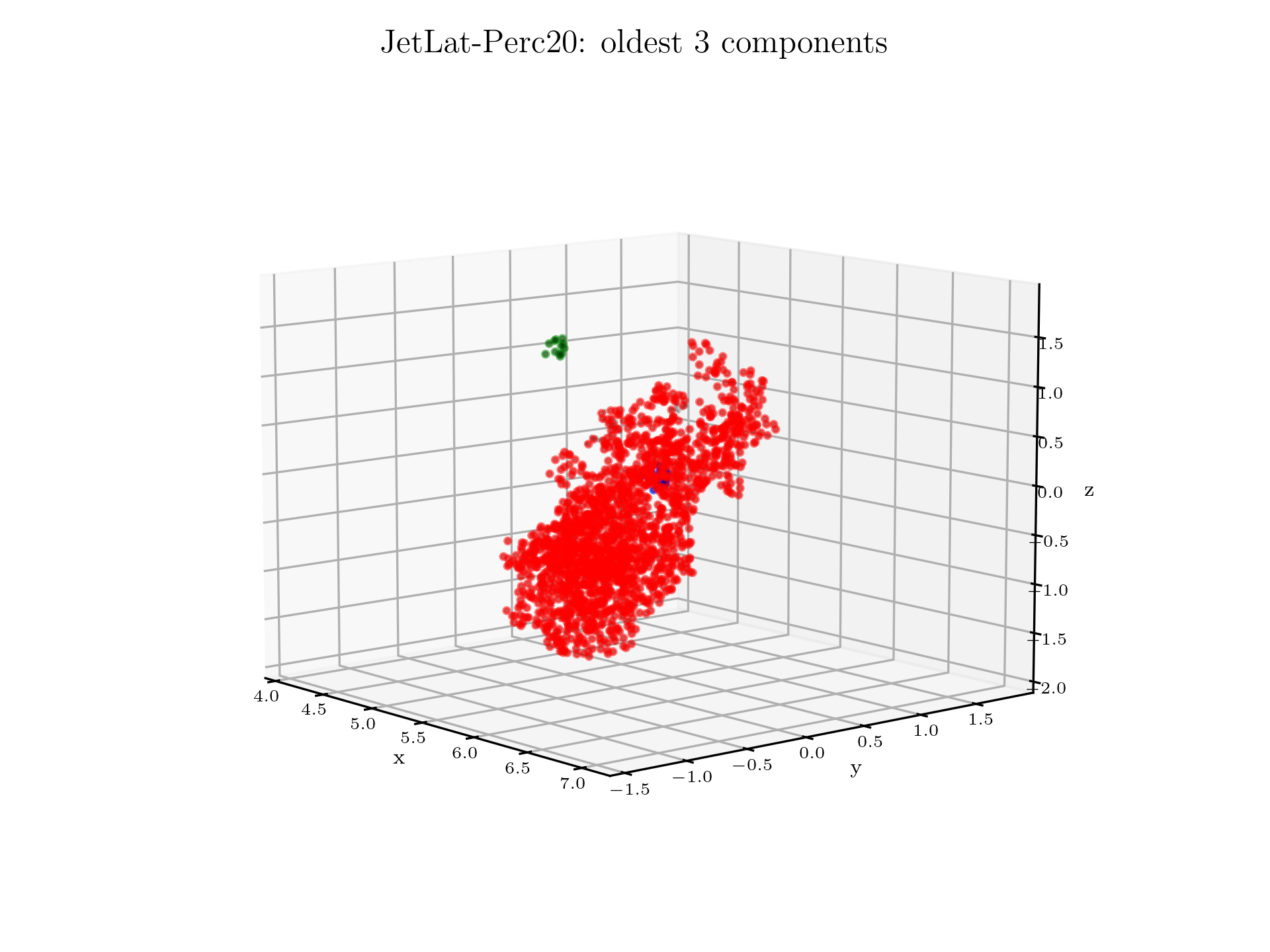}
    \caption{$C_{40}$ top 20\%}
    \label{fig:Jet_C40_20}
\end{subfigure}\hfill
\begin{subfigure}[t]{0.25\textwidth}
    \centering
    \includegraphics[width=\linewidth]{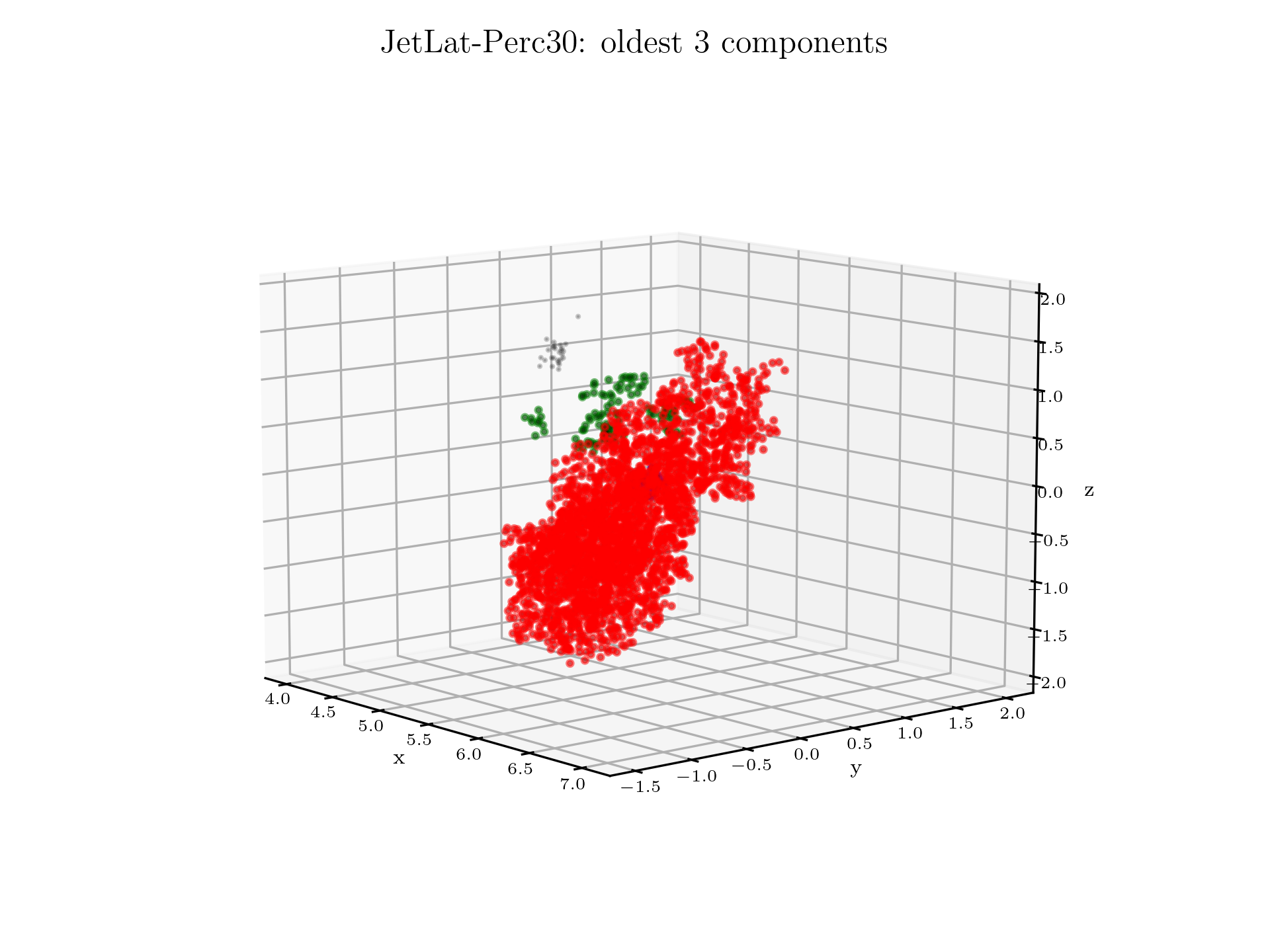}
    \caption{$C_{40}$ top 30\%}
    \label{fig:Jet_C40_30}
\end{subfigure}\hfill
\begin{subfigure}[t]{0.25\textwidth}
    \centering
    \includegraphics[width=\linewidth]{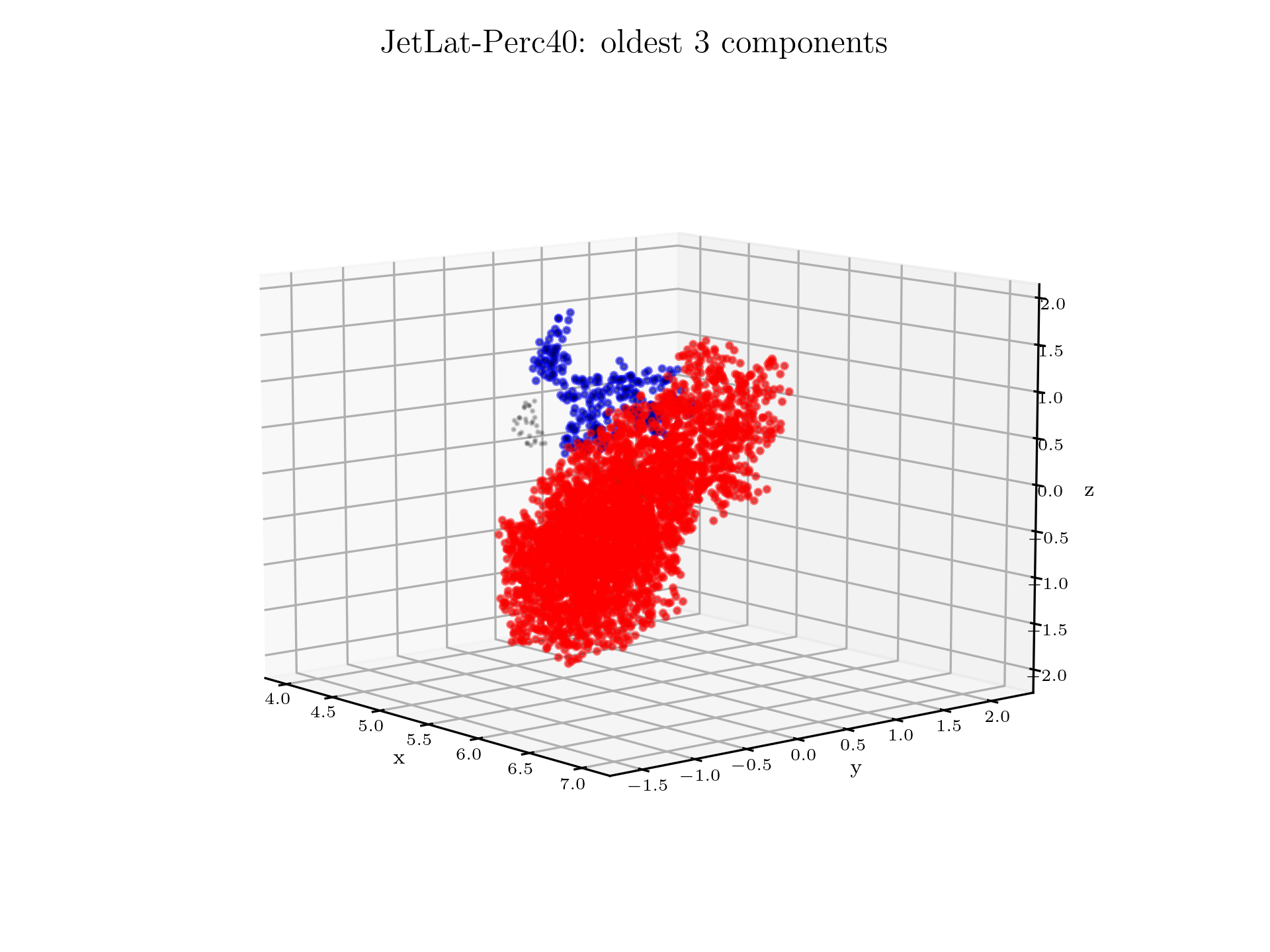}
    \caption{$C_{40}$ top 40\%}
    \label{fig:Jet_C40_40}
\end{subfigure}

\vspace{0.6em}

% ===================== Row 2: C75 =====================
\begin{subfigure}[t]{0.25\textwidth}
    \centering
    \includegraphics[width=\linewidth]{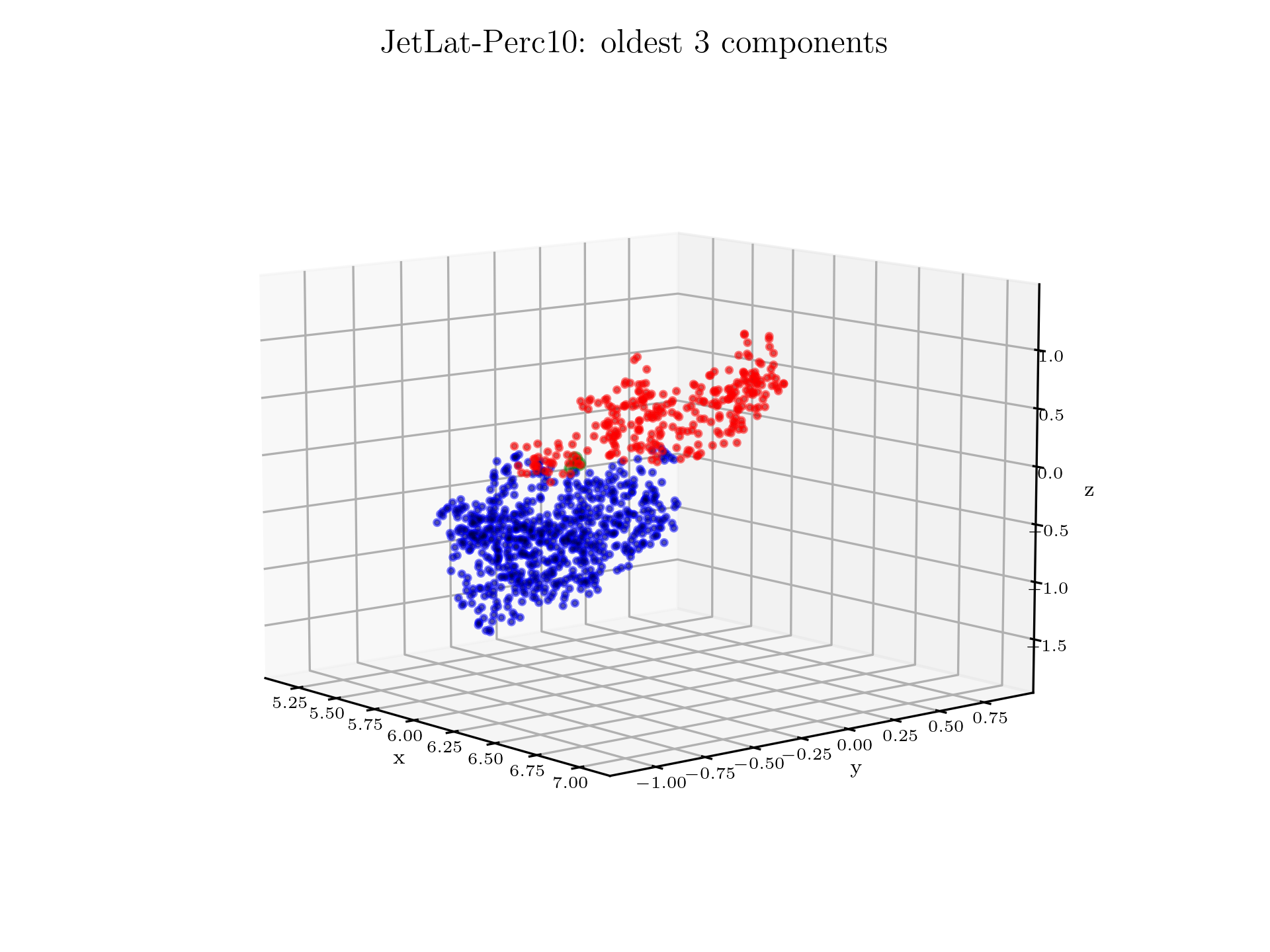}
    \caption{$C_{75}$ top 10\%}
    \label{fig:Jet_C75_10}
\end{subfigure}\hfill
\begin{subfigure}[t]{0.25\textwidth}
    \centering
    \includegraphics[width=\linewidth]{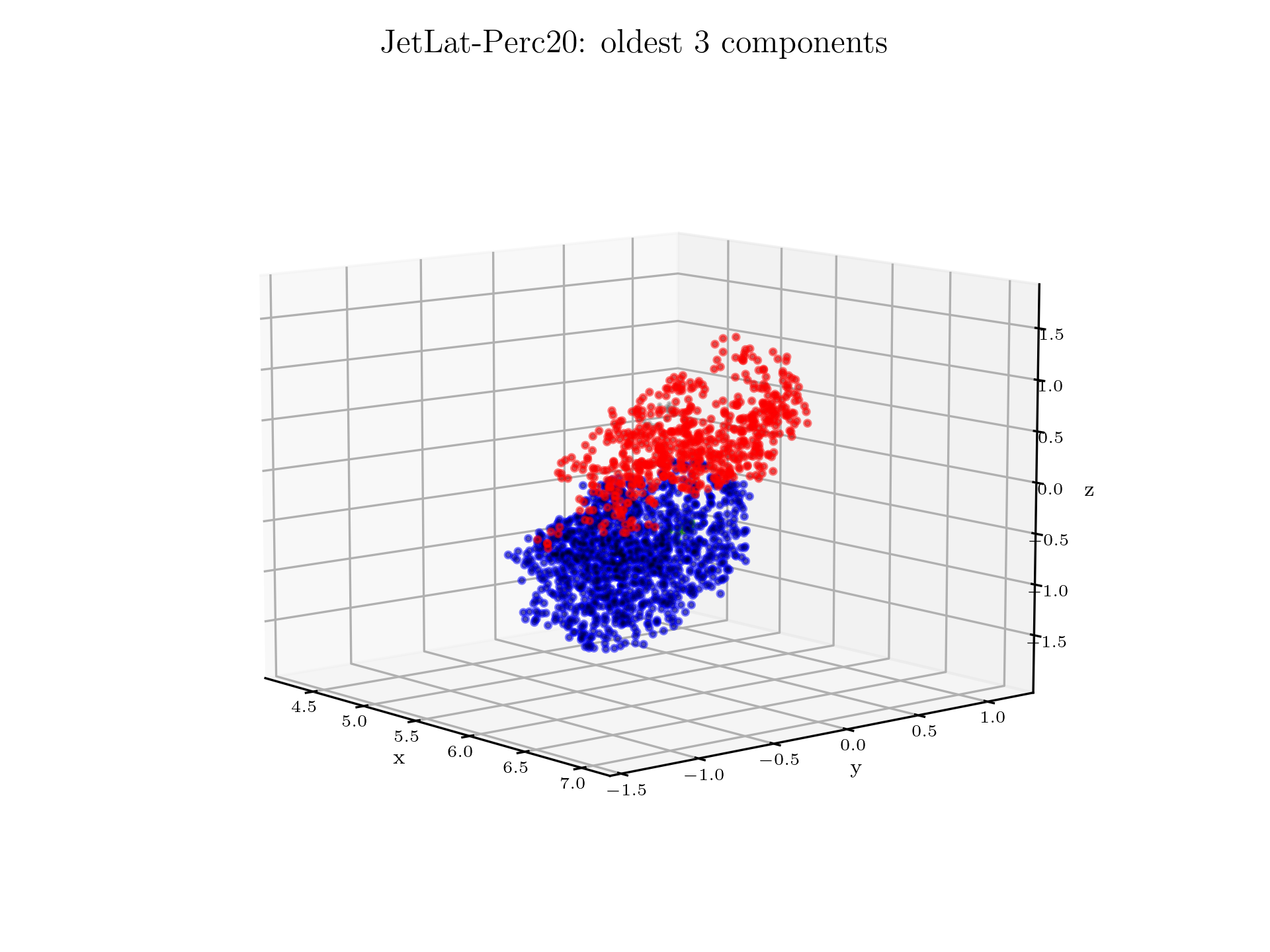}
    \caption{$C_{75}$ top 20\%}
    \label{fig:Jet_C75_20}
\end{subfigure}\hfill
\begin{subfigure}[t]{0.25\textwidth}
    \centering
    \includegraphics[width=\linewidth]{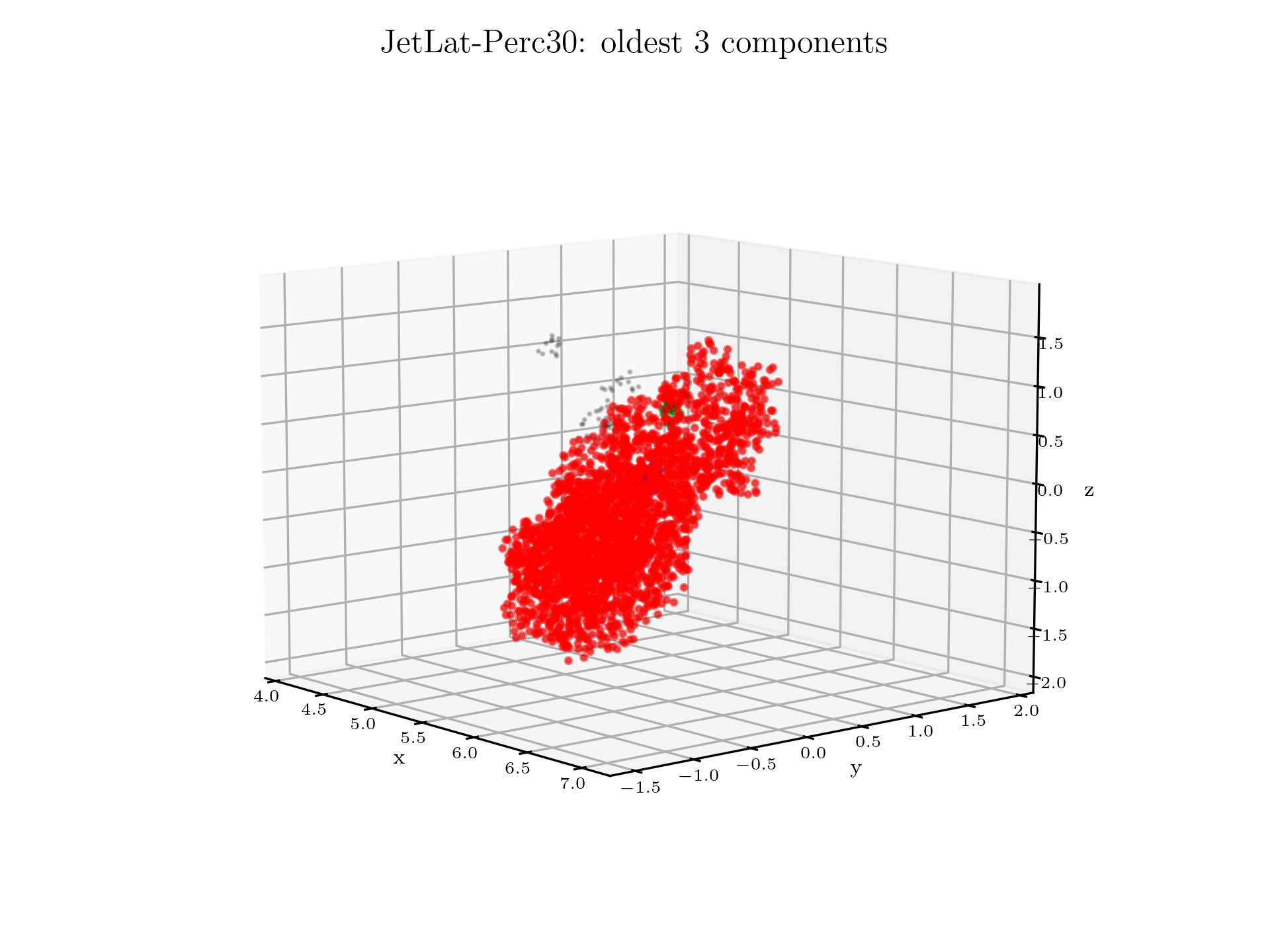}
    \caption{$C_{75}$ top 30\%}
    \label{fig:Jet_C75_30}
\end{subfigure}\hfill
\begin{subfigure}[t]{0.25\textwidth}
    \centering
    \includegraphics[width=\linewidth]{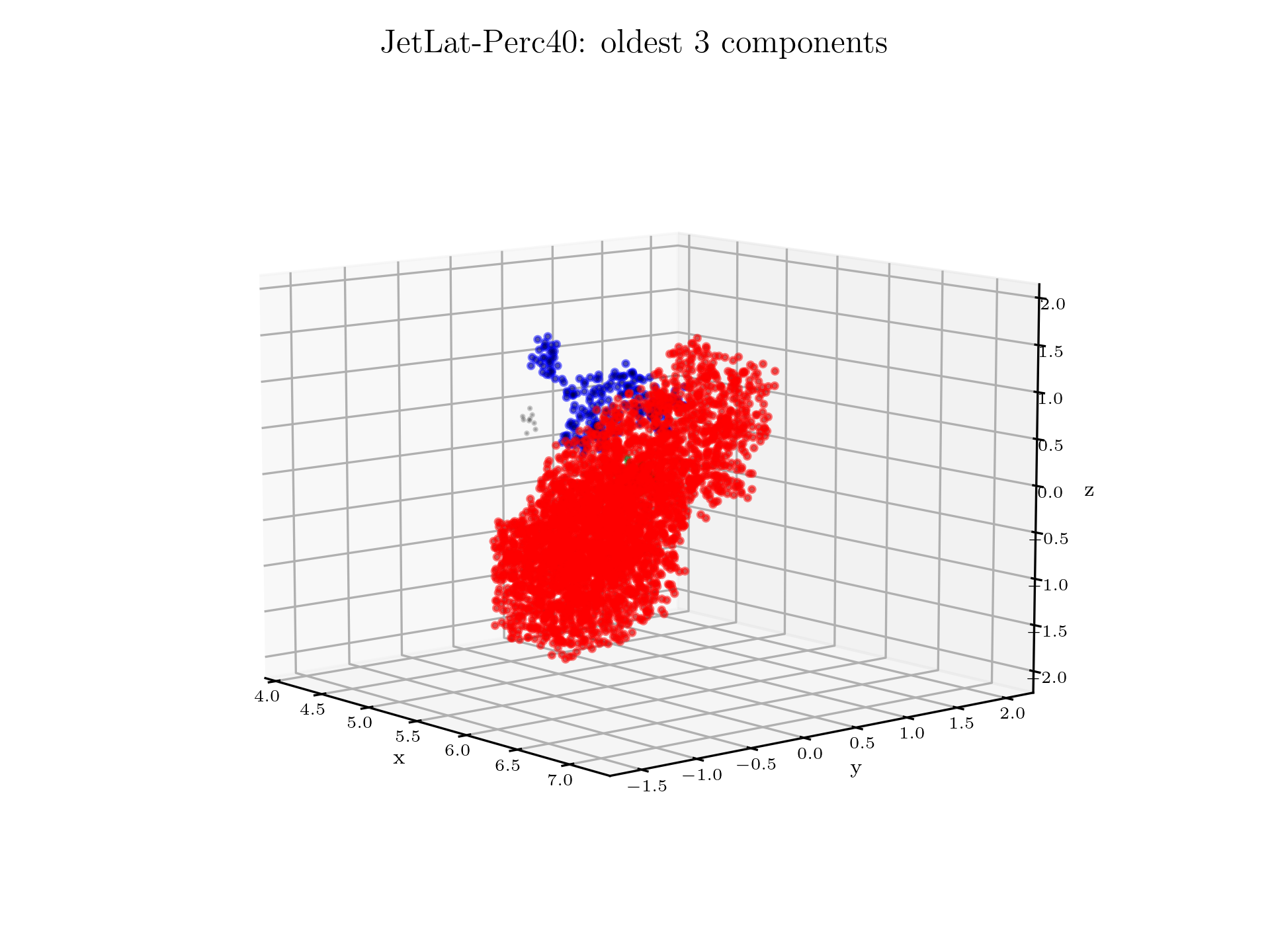}
    \caption{$C_{75}$ top 40\%}
    \label{fig:Jet_C75_40}
\end{subfigure}

\vspace{0.6em}

% ===================== Row 3: C100 =====================
\begin{subfigure}[t]{0.25\textwidth}
    \centering
    \includegraphics[width=\linewidth]{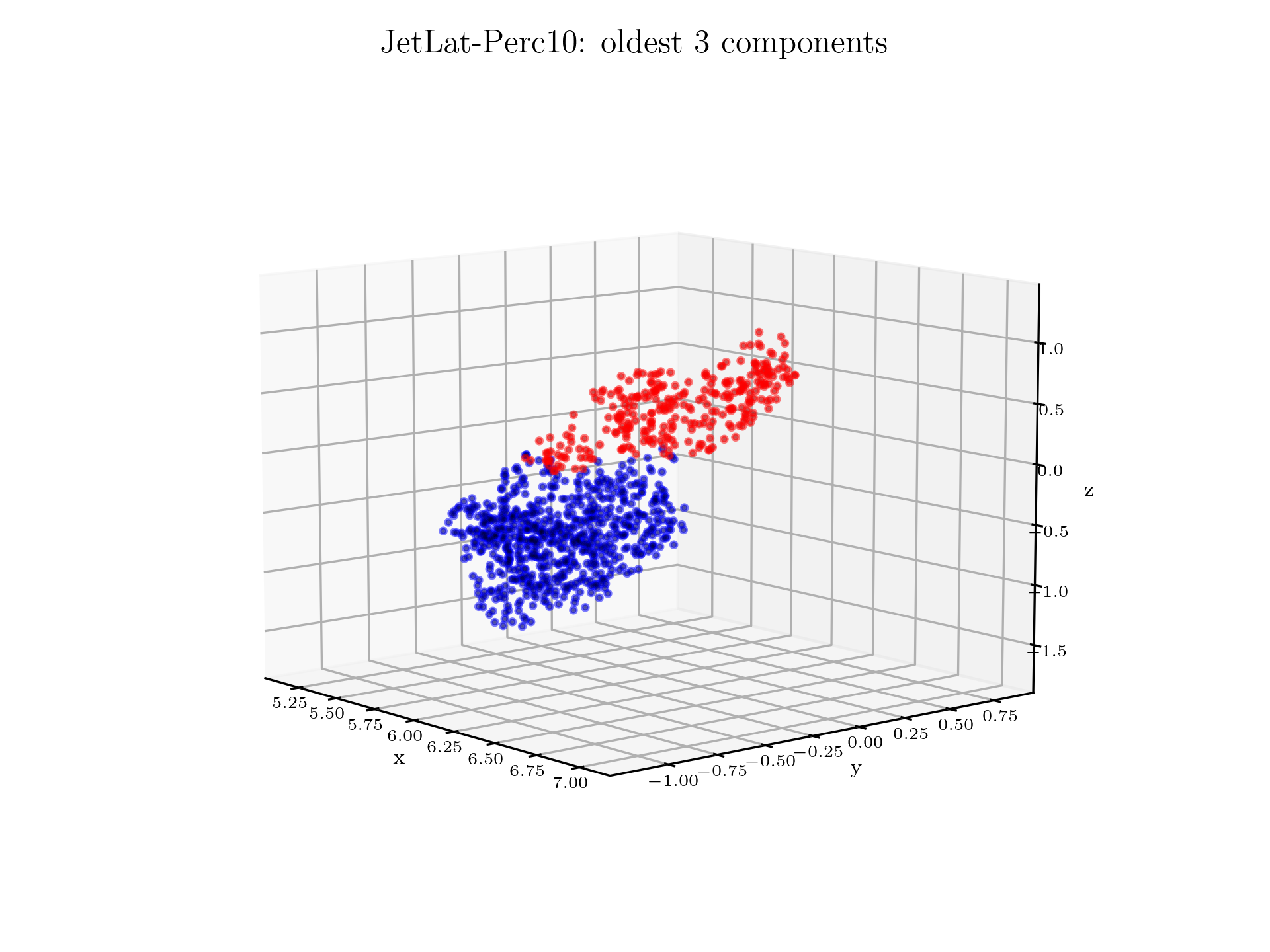}
    \caption{$C_{100}$ top 10\%}
    \label{fig:Jet_C100_10}
\end{subfigure}\hfill
\begin{subfigure}[t]{0.25\textwidth}
    \centering
    \includegraphics[width=\linewidth]{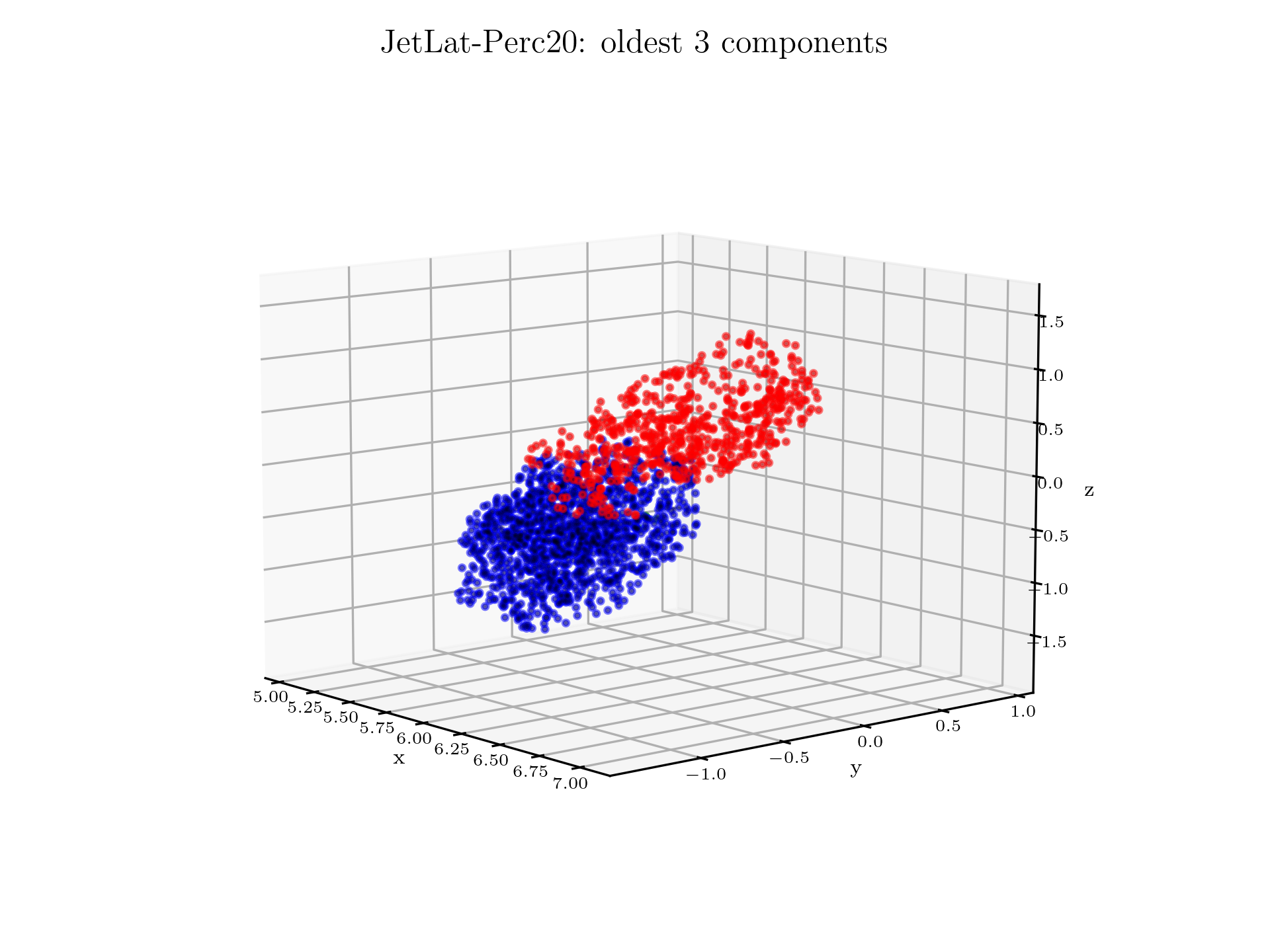}
    \caption{$C_{100}$ top 20\%}
    \label{fig:Jet_C100_20}
\end{subfigure}\hfill
\begin{subfigure}[t]{0.25\textwidth}
    \centering
    \includegraphics[width=\linewidth]{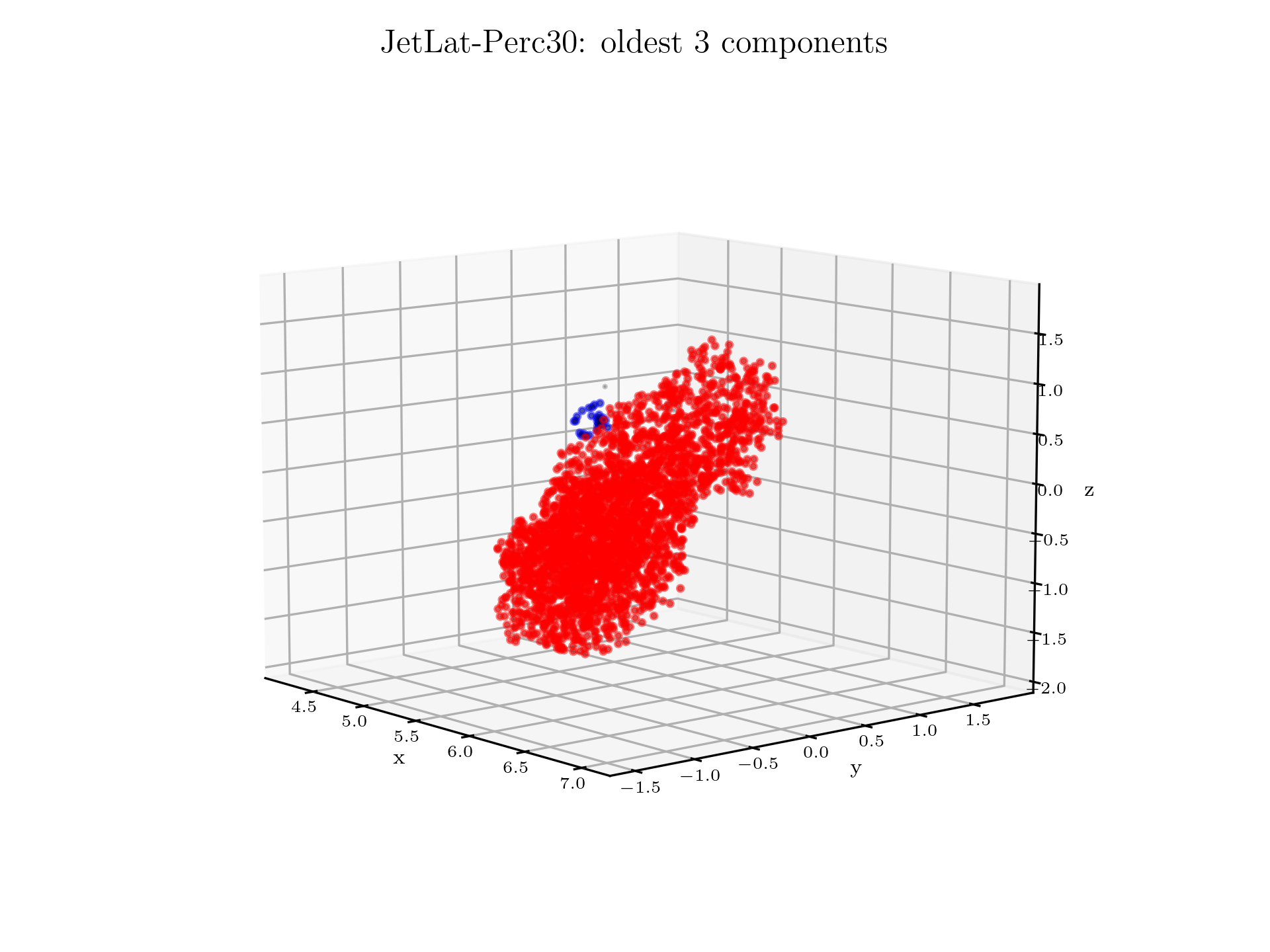}
    \caption{$C_{100}$ top 30\%}
    \label{fig:Jet_C100_30}
\end{subfigure}\hfill
\begin{subfigure}[t]{0.25\textwidth}
    \centering
    \includegraphics[width=\linewidth]{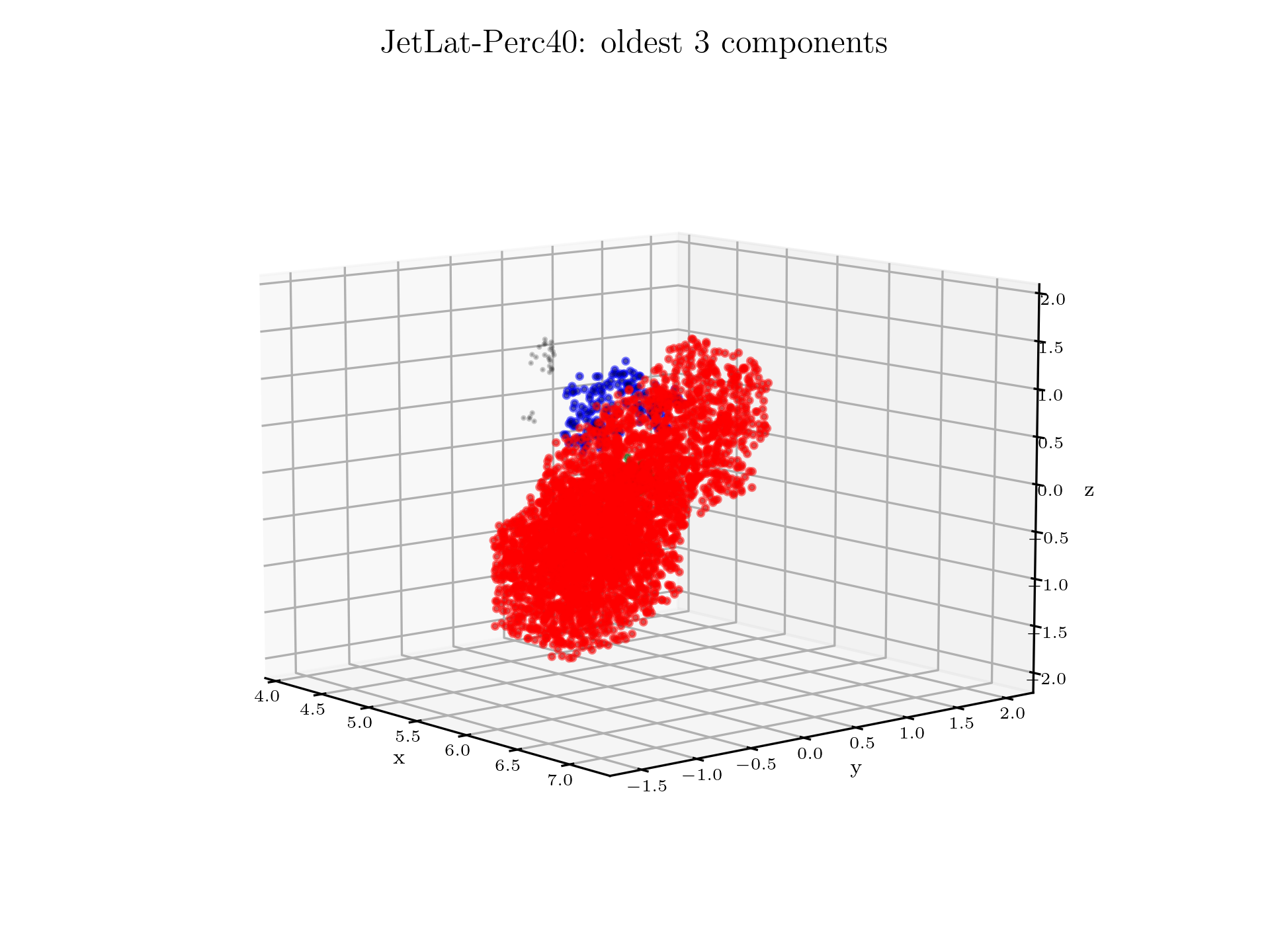}
    \caption{$C_{100}$ top 40\%}
    \label{fig:Jet_C100_40}
\end{subfigure}

\vspace{0.6em}
% ===================== Row 4: KDE =====================
\begin{subfigure}[t]{0.25\textwidth}
    \centering
    \includegraphics[width=\linewidth]{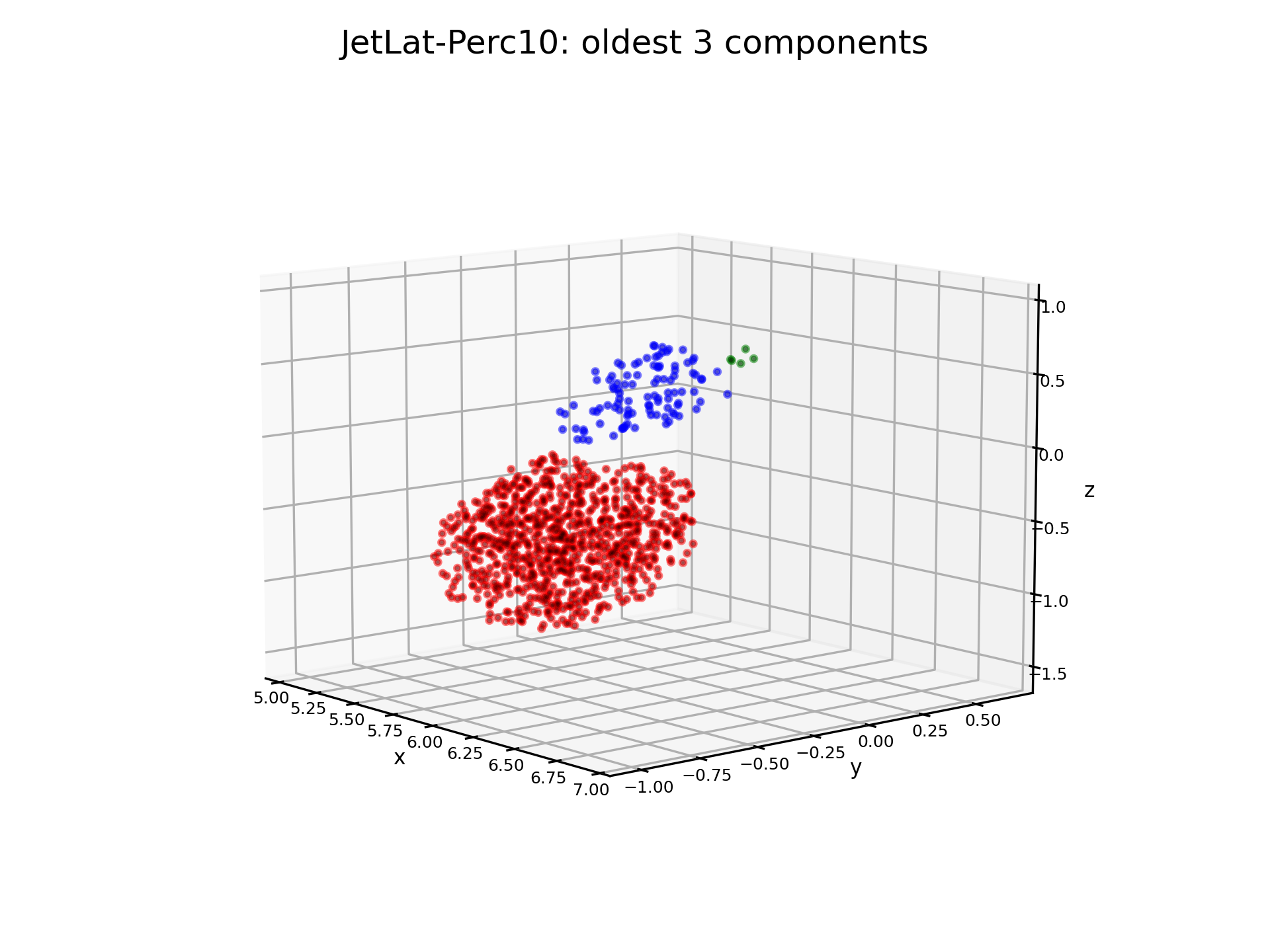}
    \caption{KDE top 10\%}
    \label{fig:Jet_KDE_10}
\end{subfigure}\hfill
\begin{subfigure}[t]{0.25\textwidth}
    \centering
    \includegraphics[width=\linewidth]{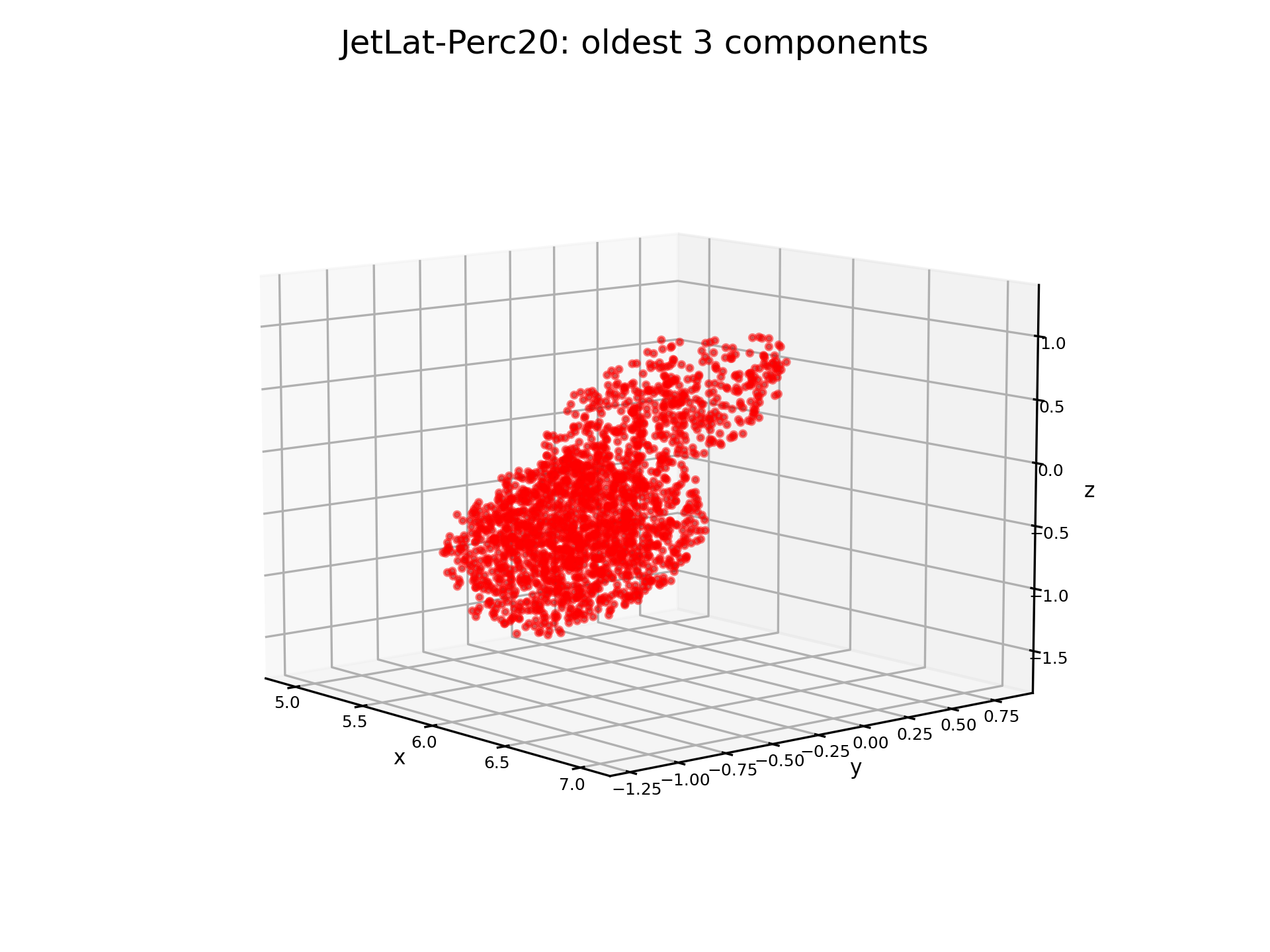}
    \caption{KDE top 20\%}
    \label{fig:Jet_KDE_20}
\end{subfigure}\hfill
\begin{subfigure}[t]{0.25\textwidth}
    \centering
    \includegraphics[width=\linewidth]{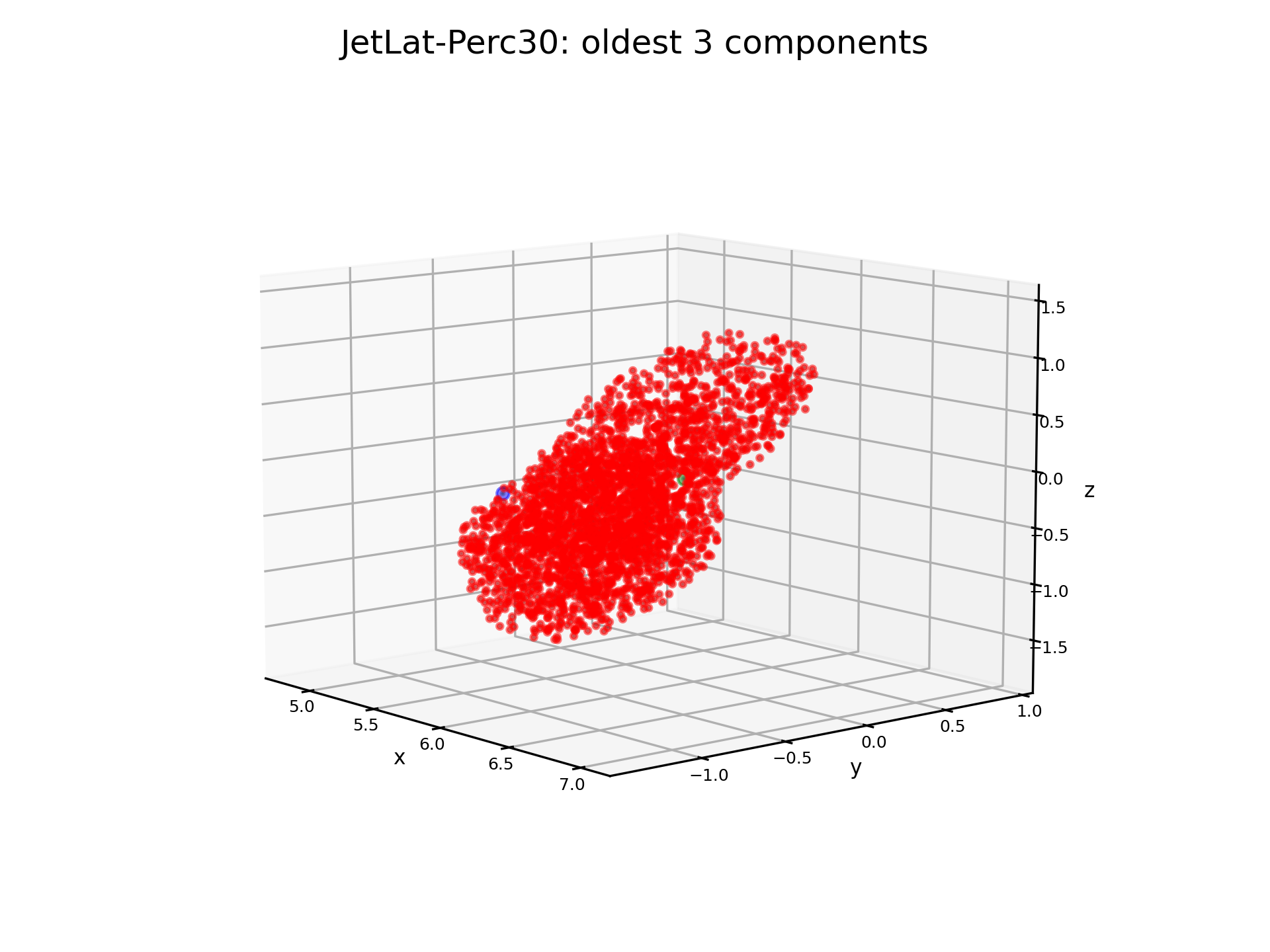}
    \caption{KDE top 30\%} 
    \label{fig:Jet_KDE_30}
\end{subfigure}\hfill
\begin{subfigure}[t]{0.25\textwidth}
    \centering
    \includegraphics[width=\linewidth]{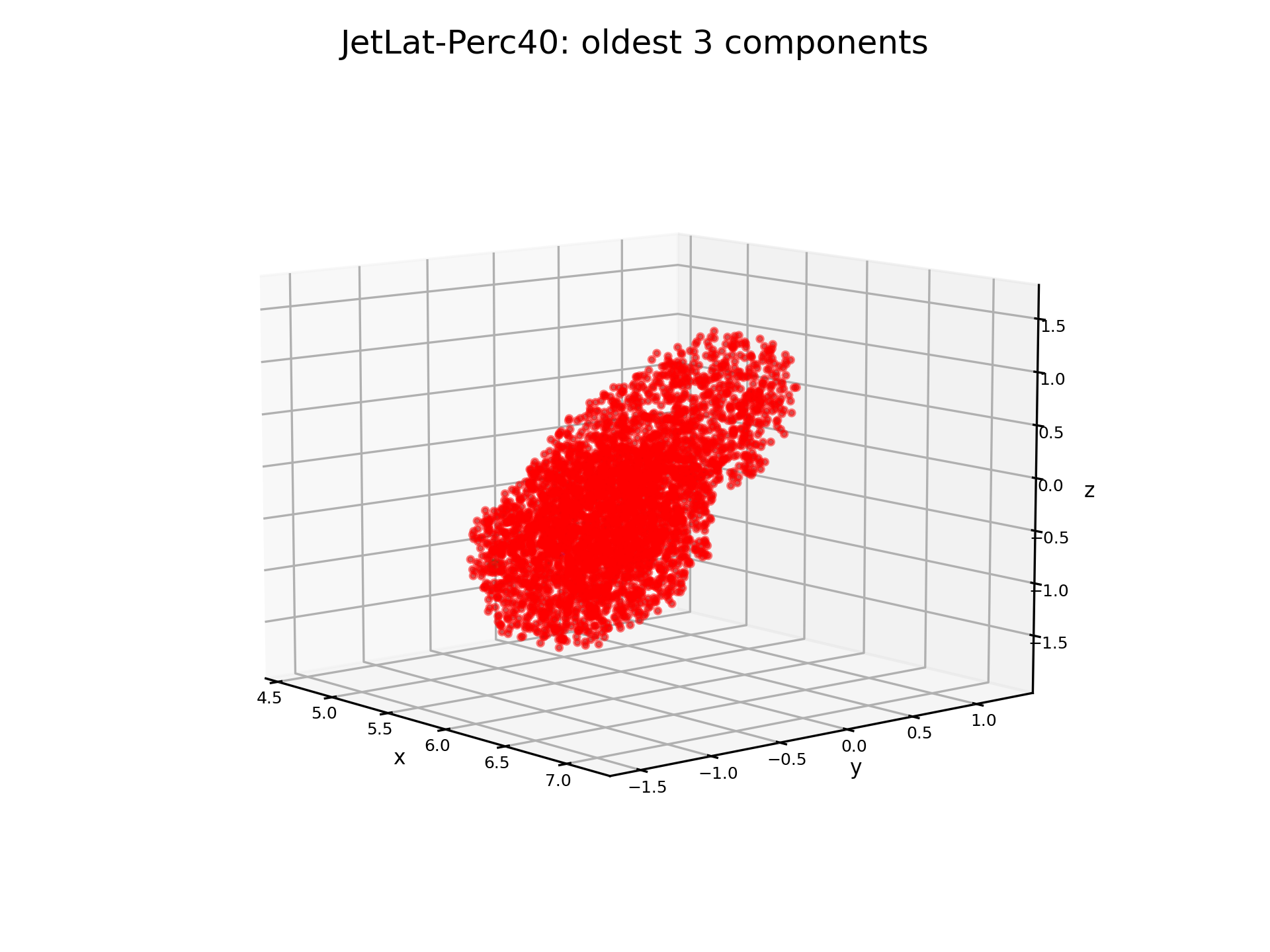}
    \caption{KDE top 40\%}
    \label{fig:Jet_KDE_40}
\end{subfigure}

\caption{
Evolution of connected components in the JetLat dataset across different bifiltration methods and percentile thresholds.
Rows correspond to filtering functions ($C_{40}$, $C_{75}$, $C_{100}$, and Gaussian KDE), while columns show increasing percentile thresholds (10\%, 20\%, 30\%, and 40\%).
Moving from left to right, connected components emerge or merge as the filtration progresses. Colors are assigned by \texttt{persloop} and indicate relative persistence of connected components, with red denoting the most persistent (longest-lived) component, followed by blue, green, and gray.
}

\label{fig:JetLat_bifiltration}
\end{figure}

Two main patterns emerge from the component structure. For $30 \leq k \leq 45$ at $p=10\%$, the three regimes are represented by three connected components, with the central regime being the most prominent and the southern regime the weakest. As $p$ increases to $20\%$, the structure changes: the central and northern regimes merge into a single large component, while the southern regime becomes fragmented into two smaller components (see Figs.~\ref{fig:Jet_C40_10}–\ref{fig:Jet_C40_40}). For larger values of $k$, a similar structure persists. The central and northern regimes are identified at low centrality thresholds and merge as the threshold increases, while the southern regime consistently appears as two components (Figs.~\ref{fig:Jet_C75_10}–\ref{fig:Jet_C75_40}).

\FloatBarrier

%====================================================
\section{Discussion, Strengths and Limitations, Future Work}\label{sec:discussion}
%====================================================
Building on the density-based topological bifiltration framework proposed by \citet{Strommen2023}, this study introduces a centrality-based extension for the topological detection of weather regimes. The original framework relies on density estimation (primarily Gaussian KDE), which can miss thin structures and loosely connected clusters. 

Our approach improves the detection of topological features in the considered datasets, in particular by better capturing thin loop structures in the CdV dataset and identifying a component that is consistent with the weakly separated southern regime in the JetLat dataset.

In addition, the centrality-based bifiltration exhibits stable regime structure over a wide range of locality parameters, suggesting robustness of the detected structures with respect to the choice of k within the tested range.

From a computational point of view, KDE and the centrality-based method differ in how they scale and represent local structure. The standard KDE implementation used in this study (and in S23) computes Gaussian kernel contributions between all pairs of data points to estimate the local density. This direct approach provides an exact and smooth estimate of the probability density but comes at a high computational cost that scales quadratically with the number of samples, $\mathcal{O}(n^2)$ \citet{silverman1986density}, which becomes restrictive for large or high-dimensional datasets.

The centrality-based method, by contrast, relies on nearest-neighbor searches rather than pairwise kernel evaluations. Efficient nearest-neighbor search requires $\mathcal{O}(n \log n)$ time for construction and supports queries in $\mathcal{O}(\log n + k)$ time per point in low to moderate dimensions \citep{bentley1975,friedman1977,omohundro1989}. Since computing the empirical dtm at each point involves averaging the distances to its $k$ nearest neighbors, evaluating the dtm for all $n$ sample points yields a total complexity of $\mathcal{O}(n \log n + nk)$, which is more efficient, particularly when $k \ll n$.

A limitation of the centrality-based approach is that, although it detects a component associated with the southern jet regime, this detection remains marginal, with relatively few points in the corresponding components. This should be taken into account when interpreting the results.

Another limitation concerns the empirical choice of the parameter $k$ for noisy datasets. Small values can make the method sensitive to noise and lead to spurious features, while large values may smooth out weak or short-lived regimes. Developing a data-driven strategy for selecting $k$ would improve the applicability of the method.

A natural future direction is to incorporate temporal information into the framework. At present, the method operates on point clouds in phase space without explicitly accounting for time ordering, which means that some of the extracted topological features may primarily reflect geometric organization rather than dynamical evolution. Extending the approach to include trajectory-based structure or temporal constraints could help isolate features that are dynamically relevant, improve the interpretation of regime transitions, and enhance the utility of the resulting representations for forecasting and statistical modeling.

%====================================================
% Declarations / Appendices (optional)
\section*{Acknowledgements}
%====================================================

This work was developed in part during the 2024 Mathematics Research Communities (MRC) program on Topological Data Analysis, supported by the National Science Foundation. We thank the MRC organizers and participants for fostering a collaborative and stimulating research environment, and we are especially grateful to Joshua Dorrington, Robyn Brooks, and Kristian Strømmen for valuable discussions and insights. We also thank the anonymous reviewers for their constructive comments and suggestions, which significantly improved the clarity, presentation, and overall quality of this work.

\section*{\bfseries Declarations}

\begin{itemize}

\item \textbf{Funding:}  
This material is based upon work supported by the National Science Foundation under Grant Number DMS-1916439.

\item \textbf{Conflict of interest:}  
The author declares that there is no conflict of interest.

\item \textbf{Ethics approval and consent to participate:}  
Not applicable.

\item \textbf{Consent for publication:}  
Not applicable.

\item \textbf{Data availability:}  
The datasets analyzed in this study are publicly available or derived from previously published sources, as described in the manuscript.

\item \textbf{Materials availability:}  
Not applicable.

\item \textbf{Code availability:}  
The analysis code used in this study is publicly available. The original implementation introduced by S23 was forked and extended to implement the centrality-based bifiltration and related analyses. The modified codebase is available at: \url{https://github.com/Soheilp86/BifiltPH}.

\item \textbf{Author contribution:}  
Soheil Anbouhi developed the methodology, performed the analysis, and wrote the manuscript.

\end{itemize}

\appendix
%==========================%L96%==========================%
\section*{Appendix: Results for the Lorenz--96 System}\label{appen:l96}
 %====================================================

\begin{figure}[ht]
    \centering
    \includegraphics[width=0.4\linewidth]{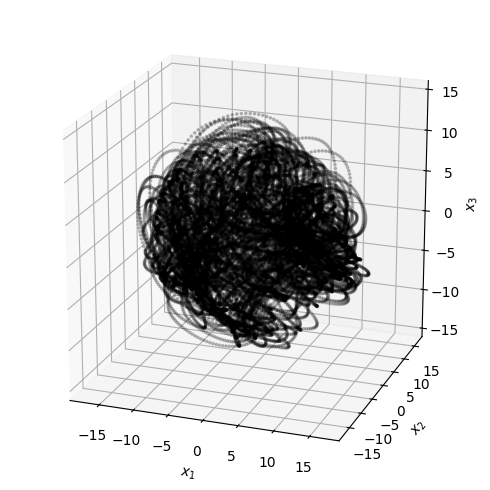}
    \caption{Lorenz--96 attractor (40{,}000 data points) projected onto its first three EOFs.}
    \label{fig:dataset-L96}
\end{figure}

The Lorenz–96 model \citet{lorenz1996predictability} is a simple system designed to mimic chaotic atmospheric dynamics. In the two-scale version used in \citep{lorenz2006regimes, christensen2015simulating, strommen2023topological} and in this study, the model consists of two coupled sets of variables: eight large-scale, low-frequency variables \(X_k\) \((k = 1, \dots, 8)\), representing dominant modes of variability, and 32 small-scale, high-frequency variables \(Y_{j,k}\) \((j = 1, \dots, 32)\), representing fast, small-scale processes. When the model trajectory is projected onto its leading empirical orthogonal functions (EOFs), the attractor exhibits two regimes: Regime~A that corresponds to a high-density region of phase space, and Regime~B that corresponds to a low-density region separating these dense clusters \citep{lorenz2006regimes, christensen2015simulating, strommen2023topological}. Fig.~\ref{fig:dataset-L96} shows the Lorenz--96 attractor in EOF space.

\begin{figure}[H]
  \centering
  % ===================== Row 1 =====================
  \begin{subfigure}[t]{0.3\textwidth}
    \centering
    \includegraphics[width=\linewidth]{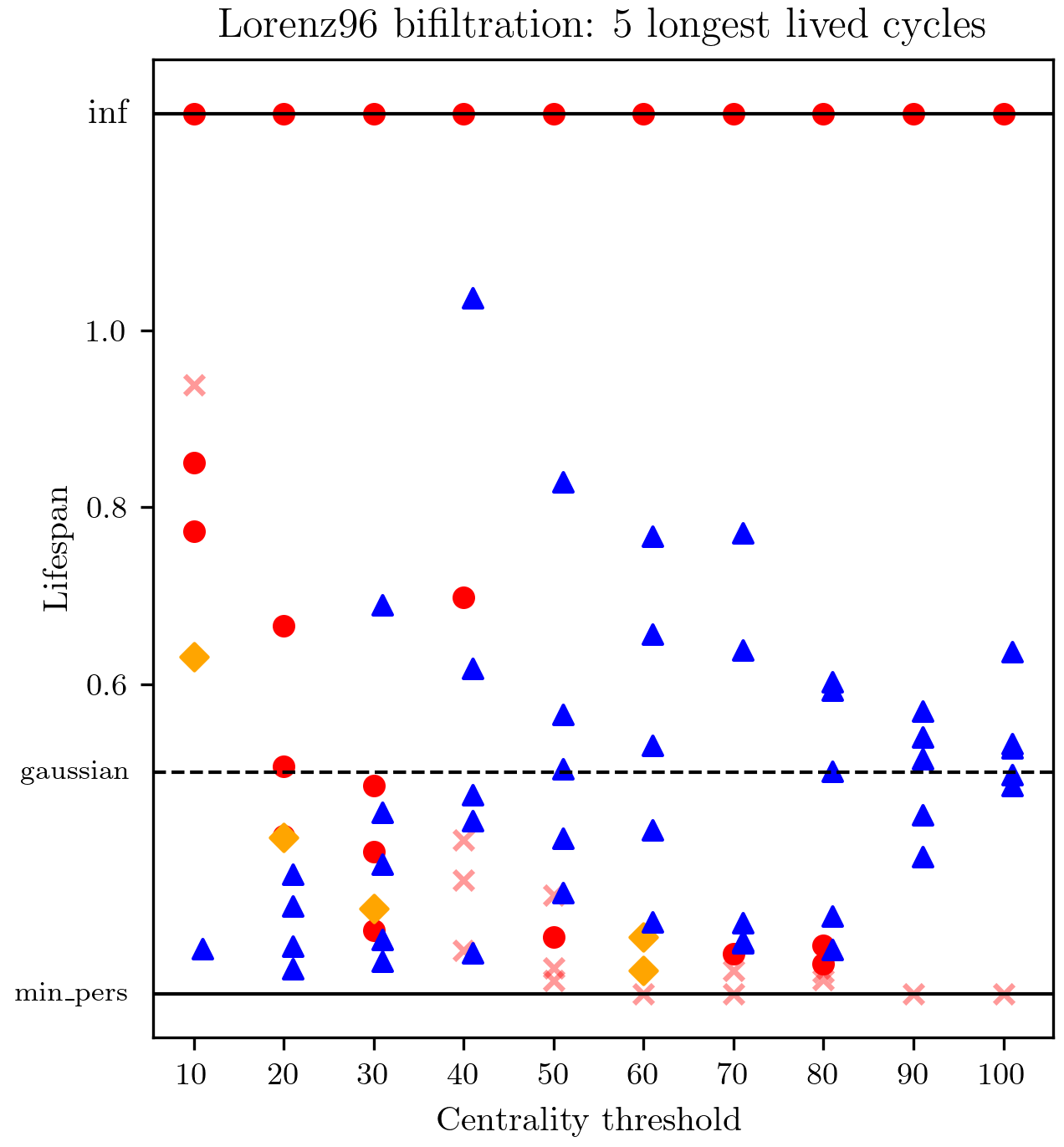}
    \caption{$C_{75}$-radius bifiltration diagram.}
    \label{fig:result_L96_C75}
  \end{subfigure}
  \hspace{0.5cm}
  \begin{subfigure}[t]{0.3\textwidth}
    \centering
    \includegraphics[width=\linewidth]{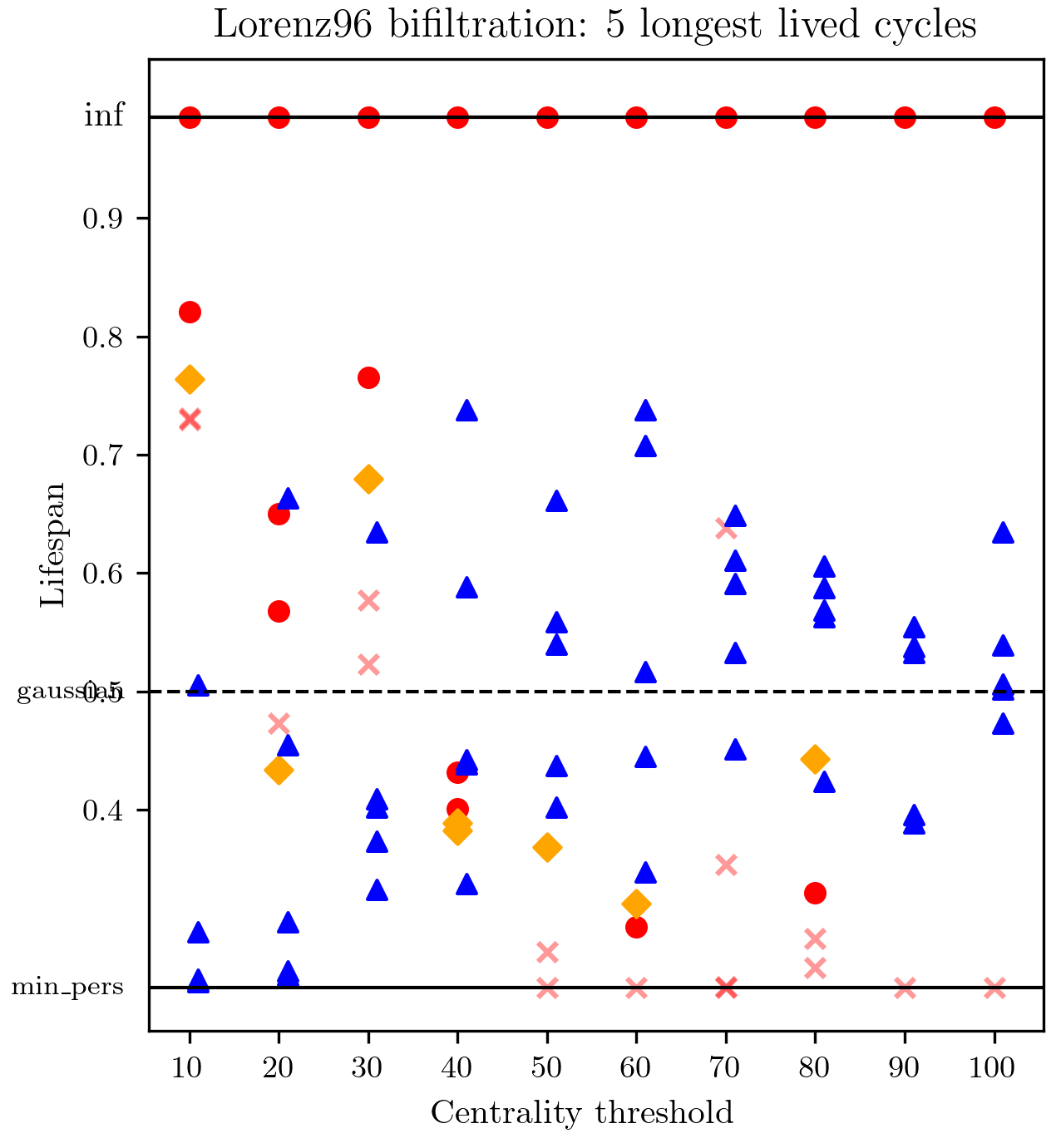}
    \caption{$C_{100}$-radius bifiltration diagram.}
    \label{fig:result_L96_C100}
  \end{subfigure}

  \vspace{0.5em}

  % ===================== Row 2 =====================
  \begin{subfigure}[t]{0.3\textwidth}
    \centering
    \includegraphics[width=\linewidth]{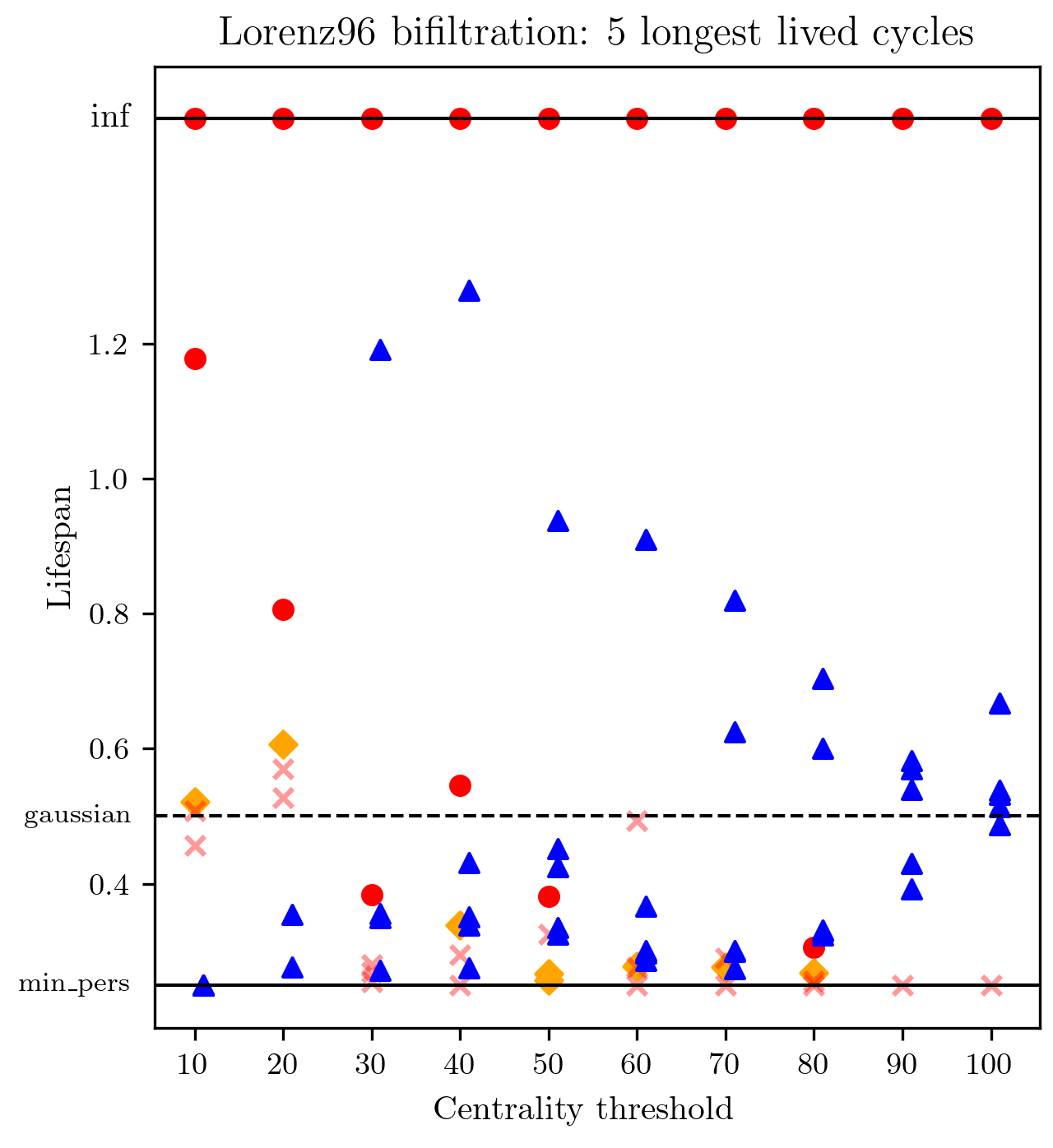}
    \caption{$C_{150}$-radius bifiltration diagram.}
    \label{fig:L96_C150}
  \end{subfigure}
  \hspace{0.5cm}
  \begin{subfigure}[t]{0.3\textwidth}
    \centering
    \includegraphics[width=\linewidth]{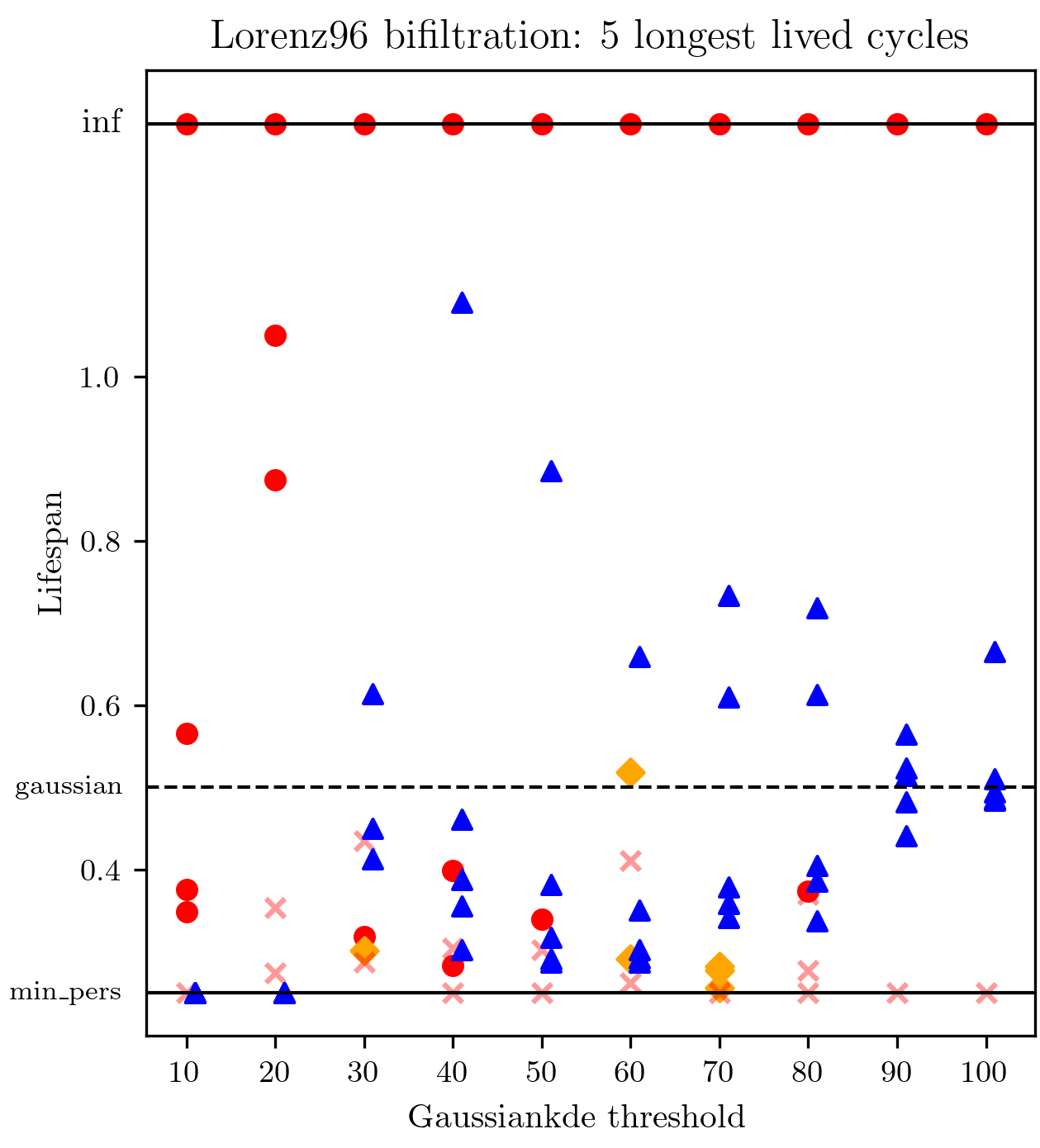}
    \caption{
        Gaussian KDE-radius bifiltration diagram}
    \label{fig:result_L96_KDE}
  \end{subfigure}

  \caption{
  Lorenz--96: comparison of centrality-based and KDE-based radius bifiltrations.
  Increasing the locality scale in the centrality functions ($C_{75}$, $C_{100}$, $C_{150}$) progressively suppresses fine-scale loop structures, leading to behavior that approaches the KDE-based bifiltration.
  }
  \label{fig:L96_result}
\end{figure}

Both KDE and centrality-based methods successfully identify the dominant regime-related structures in phase space. As in the Lorenz--63 and CdV cases, the centrality-based method with smaller \(k\) values captures these structures at lower centrality thresholds. Larger \(k\) values, on the other hand, smooth out fine details and behave similarly to the Gaussian KDE-based method. In all cases, the bifiltration first detects the densest parts of the attractor as several connected components. As the threshold increases, these components merge and the looping structures emerge. This looping behavior of the system, can be also observed in Fig.~\ref{fig:dataset-L96}.

%==========================
% Declarations / Appendices (optional)
%==========================

%\begin{appendices}
%\section{Appendix Material}\label{secA1}

%========================================
% BIBLIOGRAPHY
%========================================
\bibliographystyle{plainnat}
\bibliography{bib} 

\end{document}